\documentclass[a4paper]{article}
\usepackage{pdfsync}
\usepackage[utf8]{inputenc}
\usepackage{graphicx}
\usepackage{amssymb,amsfonts,amsmath,amsthm}
\usepackage{hyperref}
\usepackage{multirow}
\usepackage{verbatim}
\usepackage[rgb,dvipsnames]{xcolor}
\usepackage{booktabs}
\usepackage[sectionbib,round]{natbib}
\usepackage{textcomp}
\usepackage{upquote}
\usepackage{bm}
\usepackage{dirtree}

\theoremstyle{plain}

\newtheorem{proposition}{Proposition}[section]
\theoremstyle{definition}

\newtheorem{example}{Example}[section]
\newtheorem{remark}{Remark}[section]

\begin{document}

\title{\textbf{Continuous models combining slacks-based measures of efficiency and super-efficiency}}


\author{V. J. Bol\'os$^1$, R. Ben\'{\i}tez$^1$, V. Coll-Serrano$^2$ \\ \\
{\small $^1$ Dpto. Matem\'aticas para la Econom\'{\i}a y la Empresa, Facultad de Econom\'{\i}a.} \\
{\small $^2$ Dpto. Econom\'{\i}a Aplicada, Facultad de Econom\'{\i}a.} \\
{\small Universidad de Valencia. Avda. Tarongers s/n, 46022 Valencia, Spain.} \\
{\small e-mail\textup{: \texttt{vicente.bolos@uv.es}, \texttt{rabesua@uv.es}, \texttt{vicente.coll@uv.es}}} \\}

\date{July 2022}

\maketitle

\begin{abstract}
In the framework of data envelopment analysis (DEA), \cite{Ton01} introduced the slacks-based measure (SBM) of efficiency, which is a nonradial model that incorporates all the slacks of the evaluated decision-making units (DMUs) into their efficiency scores, unlike classical radial efficiency models. Next, \cite{Ton02} developed the SBM super-efficiency model in order to differentiate and rank efficient DMUs, whose SBM efficiency scores are always $1$. However, as pointed out by \cite{Che13}, some interpretation problems arise when the so-called super-efficiency projections are weakly efficient, leading to an overestimation of the SBM super-efficiency score. Moreover, this overestimation is closely related to discontinuity issues when implementing SBM super-efficiency in conjunction with SBM efficiency. \cite{Che13} and \cite{Che19} treated these problems, but they did not arrive to a fully satisfactory solution. In this paper, we review these papers and propose a new complementary score, called composite SBM, that actually fixes the discontinuity problems by counteracting the overestimation of the SBM super-efficiency score. Moreover, we extend the composite SBM model to different orientations and variable returns to scale, and propose additive versions. Finally, we give examples and state some open problems.
\end{abstract}

\section{Introduction}
\label{intro}
Data envelopment analysis (DEA) is a well known nonparametric mathematical programming technique, developed by \cite{Cha78} based on a previous work of \cite{Far57}, which allows us to assess the relative efficiency of a homogeneous set of decision-making units (DMUs) which consume several inputs in order to produce a number of outputs.
In the last two decades, there have been remarkable advances in both DEA methodologies and practical applications in a wide range of fields. Although there are several bibliographic reviews available (see, for instance, \cite{Sei97,Tav02,Emr08,Coo09}), the recent bibliographic compilation by \cite{Emr18}, providing a full listing of more than 10000 DEA-related articles ranging from 1978 to late 2016, is noteworthy.

Roughly speaking, from the observed inputs and outputs and assuming no functional relationship between them, a DEA model estimates a best-practice frontier, also known as efficient frontier, with respect to which all DMUs are evaluated. In the original CCR \citep{Cha78} and BCC \citep{Ban84} DEA radial models, the inputs are proportionally reduced while maintaining the outputs unchanged or the outputs are proportionally expanded keeping constant the inputs, depending on the orientation of the model. Shortly after, the additive model was introduced by \cite{Cha82} (see also \cite{Cha85}). This model could handle both input excesses and output shortfalls simultaneously, but it could not deliver an efficiency score as the ones obtained by the CCR and BCC radial models. To address this shortcoming, \cite{Pas99} and \cite{Ton01} proposed the enhanced Russell graph measure (ERGM) and the slacks-based measure (SBM) of efficiency respectively, which are equivalent. These models incorporate all the slacks of the evaluated DMUs into their efficiency scores, unlike classical radial efficiency models.

Usually, an efficiency DEA model classifies DMUs into two groups: efficient and inefficient. Efficient DMUs always have an efficiency score equal to $1$ but, are all efficient DMUs equally efficient?  To discriminate between efficient DMUs and rank them, \cite{And93} proposed the so-called radial super-efficiency model, whose fundamental idea is to eliminate the DMU under evaluation from the reference set. \cite{Ton02} developed the SBM super-efficiency model which consists on projecting the DMU under evaluation onto the subset of the production possibility set dominated by the DMU, and then estimating the distance between the original DMU and its projection. However, SBM super-efficiency has overestimation problems when the aforementioned projections are weakly efficient, as pointed out by \cite{Che13}. Moreover, this overestimation is closely related to discontinuity issues when implementing SBM super-efficiency together with SBM efficiency (see for example \cite{Che13, Fan13, Guo17, Che19}), contrary to what happens with radial models. In words of \cite{Che13}, this discontinuity or gap between the SBM efficiency and super-efficiency scores may lead to interpretation problems because of the sensitivity to small measurement errors or noise in the data. That is, an efficient DMU may become extremely SBM inefficient upon a small increase in inputs or a small decrease in outputs (and vice versa). 

Since continuity is a very important and desirable property for DEA models \citep{Rus90, Sch03}, the joint SBM and the continuous SBM models were introduced by \cite{Che13} and \cite{Che19} respectively, in order to solve the aforementioned discontinuity problems. However, although the joint SBM model was thought to be continuous at first, we show that, in fact, it is not always continuous. On the other hand, we also show that the continuous SBM model introduced by \cite{Che19} is not always weakly monotonic and can lead to conflicting scores in some cases. Moreover, we propose a new model, called composite SBM, that solves the discontinuity problems by counteracting the overestimation of the SBM super-efficiency score and is weakly monotonic. Nevertheless, this model presents some issues, like nonlinearity or problems related to strong monotonicity.

This paper is organized as follows. In Section \ref{sec:pre} we briefly introduce some general concepts and notation. In Section \ref{sec:sbm} we review the original SBM efficiency and super-efficiency models. In Section \ref{sec:global} we review the models presented by \cite{Che13} and \cite{Che19} in order to face the aforementioned discontinuity issue.
In Section \ref{sec:ceff0} we present the composite SBM model, studying its main properties and giving some programs for computing its score. In Section \ref{sec:ext} we extend the study to different orientations, variable returns to scale, zero or negative data, and weights. Moreover, we propose an additive version of the composite SBM model. In Section \ref{sec:examples} we give some examples showing that the interpretation and discontinuity issues are fixed by the composite SBM model. Finally, in Section \ref{sec:conc} we present some concluding remarks and state some open problems. For the sake of readability, the proofs of all the results presented in this work have been placed in Appendix \ref{sec:app} at the end of the paper.

\section{Preliminaries}
\label{sec:pre}

Notation and basic concepts are taken from \cite{Coo07}.
Vectors will be denoted by lowercase bold-face letters (either roman or greek), and they will be considered as one-column matrices when necessary. The elements of a vector will be denoted by the same letter as the vector, but unbolded and with subscripts. The $0$-vector will be denoted by $\bm{0}$ and the context will determine its dimension. All definitions and results are within the framework of constant returns to scale. Variable returns to scale are discussed in Section \ref{sec:ext}.

\subsection{Definitions}

An \emph{activity} with $m$ inputs and $s$ outputs is a pair of nonnegative vectors $\left( \mathbf{x},\mathbf{y}\right), $ where $\mathbf{x}\in \mathbb{R}_+^m$ and $\mathbf{y}\in \mathbb{R}_+^s$ are the \emph{inputs} and \emph{outputs vector} respectively. In this work, we are going to suppose that all activities are strictly positive and, therefore, the set of activities is identified with $\mathbb{R}_{>0}^{m+s}$. Nevertheless, we discuss the possibility of zero or negative data in Section \ref{sec:ext}.

Given a DMU that consumes $m$ inputs and produces $s$ outputs, it has an associated activity $\left( \mathbf{x},\mathbf{y}\right)$ given by the inputs vector $\mathbf{x}=\left( x_1,\ldots x_m\right) $ and the outputs vector $\mathbf{y}=\left( y_1,\ldots y_s\right) $, where $x_i$ is the amount of the $i$th input consumed by the DMU and $y_r$ is the amount of the $r$th output produced by the DMU, $i=1,\ldots ,m$, $r=1,\ldots ,s$. Therefore, we can identify a DMU with its activity $\left( \mathbf{x},\mathbf{y}\right)$ in the same way that a point is identified with its coordinates.
It is very important to remark that, in this work, any element of $\mathbb{R}_{>0}^{m+s}$ is called ``activity'', regardless of whether it is associated with an existing DMU or not.

Let $\mathcal{D}=\left\{ \textrm{DMU}_1, \ldots ,\textrm{DMU}_n \right\} $ be a set of $n$ DMUs, all of them having $m$ inputs and $s$ outputs. The corresponding inputs vectors $\mathbf{x}_j$, $j=1,\ldots ,n$, can be arranged as the columns of the so-called $m\times n$ \emph{input data matrix} $X$. Analogously, the outputs vectors $\mathbf{y}_j$ conform the columns of the $s\times n$ \emph{output data matrix} $Y$.
The \emph{production possibility set} defined by $\mathcal{D}$ is a set of activities given by
\begin{equation}
\label{eq:p}
P =\left\{ \left( \mathbf{x},\mathbf{y} \right) \in \mathbb{R}_{>0}^{m+s}\ \ | \ \ \mathbf{x}\geq X\bm{\lambda},\quad \mathbf{y}\leq Y\bm{\lambda},\quad \bm{\lambda}\in \mathbb{R}^n_+ \right\} ,
\end{equation}
although it is also denoted by $T$ (of \emph{Technology}) in the literature. Given two activities $\left( \mathbf{x},\mathbf{y}\right) ,\left( \mathbf{x}',\mathbf{y}'\right) $, we say that $\left( \mathbf{x},\mathbf{y}\right) $ \emph{is dominated by} $\left( \mathbf{x}',\mathbf{y}'\right) $ if $\mathbf{x}'\leq \mathbf{x}$ and $\mathbf{y}'\geq \mathbf{y}$; in this case, we say that $\left( \mathbf{x},\mathbf{y}\right) $ \emph{is strictly dominated by} $\left( \mathbf{x}',\mathbf{y}'\right) $ if $\left( \mathbf{x}',\mathbf{y}'\right) \neq \left( \mathbf{x},\mathbf{y}\right) $. The relation ``to be dominated by'' defines a partial order over the set of activities and establishes when an activity outperforms another in the sense that consumes less inputs while producing more outputs. Moreover, the production possibility set \eqref{eq:p} is formed by the activities that are dominated by \emph{positive combinations} of DMUs of the form $\left( X\bm{\lambda},Y\bm{\lambda}\right) $ with $\bm{\lambda}\in \mathbb{R}^n_+$, and hence, it is interpreted as the set of ``feasible activities'' defined by $\mathcal{D}$ \citep{Coo07}.

Given a real-valued function $f$ defined on a set of activities $\mathcal{A}$, we say that $f$ is \emph{weakly monotonic} at $\left( \mathbf{x},\mathbf{y}\right) \in \mathcal{A}$ if for any activity $\left( \mathbf{x}',\mathbf{y}'\right) \in \mathcal{A}$ such that $\left( \mathbf{x},\mathbf{y}\right) $ is dominated by $\left( \mathbf{x}',\mathbf{y}'\right) $, we have that $f\left( \mathbf{x},\mathbf{y}\right) \leq f\left( \mathbf{x}',\mathbf{y}'\right) $. Moreover, we say that $f$ is \emph{strongly monotonic} at $\left( \mathbf{x},\mathbf{y}\right) \in \mathcal{A}$ if $f\left( \mathbf{x},\mathbf{y}\right) < f\left( \mathbf{x}',\mathbf{y}'\right) $ when $\left( \mathbf{x},\mathbf{y}\right) $ is strictly dominated by $\left( \mathbf{x}',\mathbf{y}'\right) $.
We say that $f$ is \emph{weakly monotonic} on $\mathcal{A}$ if it is weakly monotonic at each activity in $\mathcal{A}$, i.e. it is order-preserving. We say that $f$ is \emph{strongly monotonic} on $\mathcal{A}$ if it is strongly monotonic at each activity in $\mathcal{A}$.

We say that an activity or a DMU is \emph{efficient} (with respect to a given set $\mathcal{D}$ of DMUs) if it is not strictly dominated by any positive combination of DMUs in $\mathcal{D}$; otherwise, we say that it is \emph{inefficient}. This concept of ``efficiency'' is equivalent to the classic ``Pareto-efficiency'' concept and it does not depend on any efficiency model. If $P$ is the production possibility set defined by $\mathcal{D}$, then any activity out of $P$ results efficient.
The set of efficient activities in $P$ is known as the \textit{(strongly) efficient frontier} (or \emph{Pareto-Koopmans frontier}) of $P$, and we denote it by $\partial ^{\text{S}}(P)$. It is clear that $\partial ^{\text{S}}(P)$ is in the frontier of $P$, known as the \emph{weakly efficient frontier} of $P$ and denoted by $\partial ^{\text{W}}(P)$. The inefficient activities in $\partial ^{\text{W}}(P)$ are known as \emph{weakly efficient}, although in fact they are not efficient.

\subsection{Score functions and efficiency scores}
\label{sec:2.2}

\begin{figure}[tbp]
\begin{center}
\includegraphics[width=0.65\textwidth]{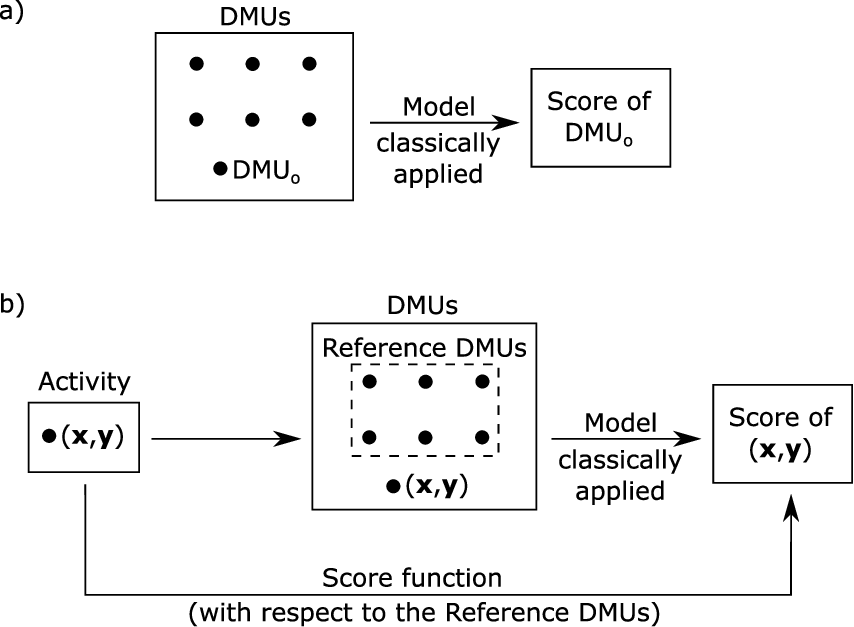}
\end{center}
\caption{a) Classically, a model is applied to a DMU (called DMU$_o$) in a given set of DMUs in order to obtain its score. b) Given a model and a set of reference DMUs, we construct a score function defined on activities. The image of an activity $\left( \mathbf{x},\mathbf{y}\right) $ is the score that the model would assign to a new hypothetical DMU with activity $\left( \mathbf{x},\mathbf{y}\right) $.}
\label{fig0}
\end{figure}

Classically, given a set of DMUs, a model is applied to one of these DMUs in order to obtain, among other things, its score (efficiency, super-efficiency, etc.). But in this work, we are going to compute scores through what we call ``score functions'': given a model and a set $\mathcal{D}$ of DMUs (which we call \emph{reference DMUs}), a \emph{score function} (with respect to $\mathcal{D}$) is a real-valued function defined on activities (i.e. from $\mathbb{R}_{>0}^{m+s}$ to $\mathbb{R}$),
such that the image of $\left( \mathbf{x},\mathbf{y}\right) $ is the score that the model would assign to a new hypothetical DMU with activity $\left( \mathbf{x},\mathbf{y}\right) $, considering $\mathcal{D}\cup \left\{ \left( \mathbf{x},\mathbf{y}\right) \right\} $ as the set of DMUs (see Figure \ref{fig0}).
There are a wide variety of models and hence, of score functions, but all of them must be at least weakly monotonic and satisfy some continuity properties, because similar activities must obtain similar scores in order to avoid sensitivity problems. Precisely, the main advantage of this methodology is that results about continuity, differentiability and monotonicity can be directly applied to score functions.

Efficiency measures (also called inefficiency measures) are the core of the DEA methodology. Restricted to inefficient activities, continuity is a property that any efficiency measure should satisfy, because discontinuities can produce serious interpretation problems \citep{Rus90, Sch03}. Moreover, monotonicity is also an important property that should be required. In this aspect, strong monotonicity is the most desirable property, but we have to note that even weak monotonicity is an elusive property for efficiency measures.
For example, \cite{And12} proved that there does not exist any weakly monotonic efficiency measure that uses a p-norm least-distance approach to the closest projections over the efficient frontier. Later, \cite{Fuk14} showed that ratio-form least-distance efficiency measures do not satisfy weak monotonicity either. \cite{And17} gave a further discussion about monotonicity of minimum distance efficiency measures. Nevertheless, any efficiency measure should be at least weakly monotonic for the sake of interpretability.

Given $\mathcal{D}$ a set of reference DMUs, an \emph{efficiency score} (with respect to $\mathcal{D}$) is a score function such that, applied to inefficient activities, represents an efficiency measure. Apart from being continuous in $P$ (the production possibility set defined by $\mathcal{D}$) and weakly monotonic, \cite{Fuk14} pointed out some other desirable properties that any efficiency score $f$ should satisfy:
\begin{enumerate}
\item $f\left( \mathbf{x},\mathbf{y}\right) =1$ if and only if $\left( \mathbf{x},\mathbf{y}\right) $ is efficient.
\item $0\leq f\left( \mathbf{x},\mathbf{y}\right) \leq 1$ for any $\left( \mathbf{x},\mathbf{y}\right) \in P$.
\end{enumerate}
Note that if $\left( \mathbf{x},\mathbf{y}\right) \notin P$, then it is efficient and hence, $f\left( \mathbf{x},\mathbf{y}\right) =1$. So,
we cannot demand global continuity, because weakly efficient activities in $\partial ^{\text{W}}(P)$ are inefficient and, according to property 1, their efficiency scores cannot be equal to $1$. Precisely, the discontinuity of efficiency scores in $\partial ^{\text{W}}(P)$ leads to discontinuity problems when implementing SBM efficiency in conjunction with SBM super-efficiency, exposed by \cite{Che13}. Note that the classical radial efficiency score is globally continuous but it does not hold property 1, since inefficient activities in $\partial ^{\text{W}}(P)$ have radial efficiency score equal to $1$ \citep{Cha62}.

\section{Original SBM models}
\label{sec:sbm}

In this Section we are going to review the SBM efficiency and super-efficiency models given by \cite{Ton01,Ton02} and define their corresponding score functions. In the rest of the paper, we are going to suppose that $\mathcal{D}$ is a set of $n$ reference DMUs, $X,Y$ are the input and output data matrices of $\mathcal{D}$ respectively, and $P$ is the production possibility set defined by $\mathcal{D}$.

We define the \emph{SBM efficiency} score (with respect to $\mathcal{D}$) of an activity $\left( \mathbf{x},\mathbf{y}\right) $ as
\begin{equation}
\label{eq:sbmact}
\def\arraystretch{1}
\begin{array}{rl}
\rho ^*\left( \mathbf{x},\mathbf{y}\right) =\min \limits_{\bm{\lambda},\mathbf{s}^-,\mathbf{s}^+} & \rho  \left( \mathbf{x},\mathbf{y},\mathbf{s}^-,\mathbf{s}^+\right) = \dfrac{1-\frac{1}{m}\sum _{i=1}^{m}s_i^- /x_{i}}{1+\frac{1}{s}\sum _{r=1}^{s}s_r^+ /y_{r}} \\
\text{s.t.} & \mathbf{x} = X\bm{\lambda}+\mathbf{s}^-, \\
& \mathbf{y} = Y\bm{\lambda}-\mathbf{s}^+, \\
& \bm{\lambda}\in \mathbb{R}^n_+,\quad \mathbf{s}^-\in \mathbb{R}^m_+,\quad\mathbf{s}^+\in \mathbb{R}^s_+,
\end{array}
\end{equation}
if $\left( \mathbf{x},\mathbf{y}\right) \in P$, and $\rho ^*\left( \mathbf{x},\mathbf{y}\right) =1$ if $\left( \mathbf{x},\mathbf{y}\right) \notin P$. The vectors $\mathbf{s}^-,\mathbf{s}^+$ are called \emph{inefficiency slack} vectors \citep{Guo17}. Considering $\mathbf{s}^{-*},\mathbf{s}^{+*}$ optimal inefficiency slack vectors, we refer to the activities of the form $\left( \mathbf{x} -\mathbf{s}^{-*}, \mathbf{y} +\mathbf{s}^{+*}\right) $ as \emph{efficient} (or \emph{optimal}) \emph{targets}. Program \eqref{eq:sbmact} is based on the original SBM efficiency model given by \cite{Ton01}, in the sense that $\rho ^*\left( \mathbf{x},\mathbf{y}\right) $ is the score that the Tone's original SBM efficiency model would assign to a new DMU with activity $\left( \mathbf{x},\mathbf{y}\right)$. 
We conclude that $\rho ^*$ is an efficiency score because it satisfies properties 1 and 2, it is weakly monotonic and $\rho ^*|_P$ is clearly continuous. With respect to monotonicity, we have the next result:
\begin{proposition}
\label{prop:rhomon}
The SBM efficiency score $\rho ^*|_P$ is strongly monotonic.
\end{proposition}

We define the \emph{SBM super-efficiency} (S-SBM) score (with respect to $\mathcal{D}$) of an activity $\left( \mathbf{x},\mathbf{y}\right) $ as
\begin{equation}
\label{eq:ssbmact}
\def\arraystretch{1}
\begin{array}{rl}
\delta ^*\left( \mathbf{x},\mathbf{y}\right) =\min \limits_{\bm{\lambda},\mathbf{t}^-,\mathbf{t}^+} & \delta \left( \mathbf{x},\mathbf{y},\mathbf{t}^-,\mathbf{t}^+\right)= \dfrac{1+\frac{1}{m}\sum _{i=1}^{m}t_i^- /x_{i}}{1-\frac{1}{s}\sum _{r=1}^{s}t_r^+ /y_{r}} \\
\text{s.t.} & \mathbf{x}+\mathbf{t}^- \geq X\bm{\lambda}, \\
& \mathbf{0} < \mathbf{y}-\mathbf{t}^+ \leq Y\bm{\lambda}, \\
& \bm{\lambda}\in \mathbb{R}^n_+,\quad \mathbf{t}^-\in \mathbb{R}^m_+,\quad\mathbf{t}^+\in \mathbb{R}^s_+,
\end{array}
\end{equation}
where $\mathbf{t}^-,\mathbf{t}^+$ are called \emph{super-efficiency slack} vectors \citep{Guo17}. Taking into account the SBM super-efficiency model given by \cite{Ton02}, we have that $\delta ^*\left( \mathbf{x},\mathbf{y}\right) $ is the score that a new DMU with activity $\left( \mathbf{x},\mathbf{y}\right) $ would have.
Note that in Tone's S-SBM program, the DMU under evaluation must be excluded from the original set of DMUs. However, in program \eqref{eq:ssbmact}, there is no need to make any exclusion. This is one of the advantages of working with score functions.

The set of activities in $P$ that are dominated by $\left( \mathbf{x},\mathbf{y}\right) $ is given by
\begin{equation}
\label{eq:barP}
\bar{P}=\left\{ \left( \bar{\mathbf{x}},\bar{\mathbf{y}}\right) \in P \ \ | \ \ \bar{\mathbf{x}}\geq \mathbf{x},\quad \bar{\mathbf{y}}\leq \mathbf{y}\right\} .
\end{equation}
Hence, the constraints of \eqref{eq:ssbmact} are equivalent to $\left( \mathbf{x}+\mathbf{t}^-,\mathbf{y}-\mathbf{t}^+\right) \in \bar{P}$ and then, according to \cite{Ton02}, $\delta $ can be interpreted as a weighted $l_1$ distance from $\left( \mathbf{x},\mathbf{y}\right) $ to activities in $\bar{P}$. Note that although technically $\delta $ is not a distance in the mathematical sense, we are going to use the term ``distance'' in the same way that it is used in \cite{Ton02}. In this case, if $\mathbf{t}^{-*},\mathbf{t}^{+*}$ are optimal super-efficiency slack vectors for \eqref{eq:ssbmact}, then $\left( \mathbf{x}+\mathbf{t}^{-*},\mathbf{y}-\mathbf{t}^{+*}\right) $ are the activities in $\bar{P}$ closest to $\left( \mathbf{x},\mathbf{y}\right) $, which are called \emph{super-efficiency projections} of $\left( \mathbf{x},\mathbf{y}\right) $. Note that, since the distance between a point and a closed set is defined by the distance between the point and the closest point in the set, $\delta ^*\left( \mathbf{x},\mathbf{y}\right) $ can also be interpreted as a distance from $\left( \mathbf{x},\mathbf{y}\right) $ to $\bar{P}$.
It is important to remark that, according to \cite{Che13}, an overestimation of the S-SBM score $\delta ^*\left( \mathbf{x},\mathbf{y}\right) $ is produced when $\left( \mathbf{x},\mathbf{y}\right) $ has  weakly efficient (and hence inefficient) super-efficiency projections. This overestimation occurs because the inefficiency of such projections is not taken into account by $\delta ^*\left( \mathbf{x},\mathbf{y}\right) $. This fact is closely related to discontinuity issues when implementing SBM efficiency in conjunction with SBM super-efficiency, as we will see in Section \ref{sec:global}.

It is clear that the S-SBM score $\delta ^*$ is continuous. With respect to monotonicity,
\cite{Ton02} proved that $\delta ^*\left( \alpha \mathbf{x}, \beta \mathbf{y}\right) \geq \delta^*\left( \mathbf{x}, \mathbf{y}\right) $ for any activity $\left( \mathbf{x}, \mathbf{y}\right) $, $\alpha \leq 1$ and $\beta \geq 1$.
Next, we present a more general result.

\begin{proposition}
\label{prop:22}
The S-SBM score $\delta ^*$ is weakly monotonic.
\end{proposition}

\begin{remark}[Strong monotonicity]
\label{rem:ssbmnm}
The S-SBM score $\delta ^*$ is constantly equal to $1$ for activities in $P$ and hence, it is obvious that it is not strongly monotonic for inefficient activities. However, it is important to note that $\delta ^*$ is not strongly monotonic for efficient activities either, because it does not take into account the inefficiency of super-efficiency projections. For example, let us consider $\mathcal{D}=\left\{ \text{D1},\text{D2},\text{D3}\right\} $ a set of reference DMUs with two inputs and one normalized output, where $\text{D1}=\left( \left( 10,40\right) ,1\right) $, $\text{D2}=\left( \left( 15,25\right) ,1\right) $ and $\text{D3}=\left( \left( 30,20\right) ,1\right) $. As we will see in Example \ref{ex81}, activities of the form  $\left( \left( 30+c,x_2\right) ,1\right) $ with $0<x_2<20$ have the same S-SBM score $\delta ^*$ (with respect to $\mathcal{D}$) when $c\geq 0$ varies.
\end{remark}

\section{Global continuous SBM models}
\label{sec:global}

From the SBM efficiency and super-efficiency models defined by \cite{Ton01} and \cite{Ton02} respectively, it is possible to construct a global SBM score function defined on the whole $\mathbb{R}_{>0}^{m+s}$ such that it coincides with $\rho^*$ for activities in the corresponding production possibility set $P$ and, on the other hand, it coincides with $\delta^*$ for the rest of activities (see \cite{Fan13,Guo17}). Nevertheless, as it was firstly showed by \cite{Che13}, this score is not continuous in the weakly efficient frontier $\partial ^{\text{W}}(P)$, making it hard to interpret and justify the scores in applications. As it is also pointed out by \cite{Che13}, this discontinuity issue is closely related to the overestimation of the S-SBM score produced when the super-efficiency projections are weakly efficient. So, the idea to fix the discontinuity problem is to define a new model that penalizes the inefficiency of such projections in some way.

In this section we are going to review the joint SBM (J-SBM) and the continuous SBM (CSBM) models introduced by \cite{Che13} and \cite{Che19} respectively. Both models try to solve the discontinuity problem but, unfortunately, they do not give a fully satisfactory solution for different reasons that we are going to show. Specifically, the J-SBM score is not continuous in some cases and the CSBM score, even being continuous, is not weakly monotonic in some cases.
It must be taken into account that we have changed some notation to simplify and adapt it to our study.

\subsection{The joint SBM model}

Given $\mathcal{D}$, the space of activities $\mathbb{R}_{>0}^{m+s}$ is splitted into three regions:
\begin{itemize}
\item{Region (I)} or \emph{technical inefficiency zone}: $P-\partial^{\text{S}}(P)$, i.e. inefficient activities.
\item{Region (II)}: efficient activities with all super-ef\-fi\-cien\-cy projections being efficient.
\item{Region (III)}: efficient activities with at least one inefficient super-efficien\-cy projection.
\end{itemize}
With respect to the notation of \cite{Che13}, super-efficiency projections are called ``S-SBM reference points'' and, analogously for inefficient activities, efficient targets are called ``SBM reference points''.

Given an activity $\left( \mathbf{x},\mathbf{y}\right) $, the J-SBM model assigns a score $\phi ^* \left( \mathbf{x},\mathbf{y}\right) $ that coincides with $\rho ^*\left( \mathbf{x},\mathbf{y}\right) $ if $\left( \mathbf{x},\mathbf{y}\right) $ is in Region (I), is equal to $\delta ^*\left( \mathbf{x},\mathbf{y}\right) $ if $\left( \mathbf{x},\mathbf{y}\right) $ is in Region (II) and, for activities in Region (III), the author introduces a \textit{relaxed SBM super-efficiency} model that penalizes the inefficiency of super-efficiency projections:
\begin{equation}
\label{eq:relaxingcup}
\def\arraystretch{1}
\begin{array}{rl}
\phi ^*\left( \mathbf{x},\mathbf{y}\right) =\min \limits_{\bm{\lambda},\tilde{\mathbf{s}}^-,\tilde{\mathbf{s}}^+} & \phi \left( \mathbf{x},\mathbf{y},\tilde{\mathbf{s}}^-,\tilde{\mathbf{s}}^+\right)= \dfrac{1-\frac{1}{m}\sum _{i=1}^{m}\tilde{s}_i^- /x_{i}}{1+\frac{1}{s}\sum _{r=1}^{s} \tilde{s}_r^+/y_{r}} \\
\text{s.t.} & X\bm{\lambda} = \mathbf{x}-\tilde{\mathbf{s}}^-, \\
& Y\bm{\lambda}=\mathbf{y}+\tilde{\mathbf{s}}^+, \\
& \bm{\lambda}\in \mathbb{R}^n_+,\quad \tilde{\mathbf{s}}^-\in \mathbb{R}^m \textrm{ free},\quad \tilde{\mathbf{s}}^+\in \mathbb{R}^s\textrm{ free}.
\end{array}
\end{equation}
The \emph{reference points} given by model \eqref{eq:relaxingcup} are those of the form $\left( \mathbf{x}-\tilde{\mathbf{s}}^{-*},\mathbf{y}+\tilde{\mathbf{s}}^{+*}\right) $ where $\tilde{\mathbf{s}}^{-*}$, $\tilde{\mathbf{s}}^{+*}$ are optimal. Note that $\tilde{\mathbf{s}}^-$, $\tilde{\mathbf{s}}^+$ are free slack vectors.

\cite{Che13} uses binary variables in order to express the J-SBM score in a single model, constructing a ``switch'' between the three different models from each region (see Equation (9) in \cite{Che13}). However, in the way the model is expressed, the ``switch'' does not work correctly between Region (II) and Region (III). Anyway, this mistake could be fixed by re-defining the J-SBM score piecewisely instead of using binary variables, although its computation would need several stages. Namely, in a first stage, we need to know to which region the activity belongs and then, in a second stage, we apply \eqref{eq:sbmact}, \eqref{eq:ssbmact} or \eqref{eq:relaxingcup} for activities in Region (I), (II) or (III), respectively.

Moreover, Corollary 1 in \cite{Che13}, which states that the reference points given by the J-SBM model are efficient, is not fulfilled.
This issue is treated by \cite{Lin18}, giving a counterexample and a revised model.

Finally, and most importantly, according to Theorem 5 in \cite{Che13}, the J-SBM score is supposed to connect all three regions in a continuous way. Unfortunately, this is not true in some cases as we show in Example \ref{ex81} (see Figure \ref{figmapas_zoom} (b)), and this mistake can not be fixed. The reason for this to happen is that, sometimes, the relaxed model \eqref{eq:relaxingcup} gives a reference point quite far away from the evaluated activity. The discontinuity of the J-SBM score was recognised in a corrigendum paper by \cite{Che14}, ensuring that ``it can be easily corrected by constraining the reference point to be fixated on a specific strongly Pareto-efficient point''. However, the details of this proposition were not given and, as far as we know, there is not any further paper that clarifies it.

\subsection{The continuous SBM model}

Given an efficient activity $\left( \mathbf{x},\mathbf{y}\right) $ with $\mathbf{t}^{-*},\mathbf{t}^{+*}$ optimal super-efficiency slack vectors for \eqref{eq:ssbmact}, the SBM efficiency score (with respect to $\mathcal{D}$) of its super-efficiency projection given by $\left( \mathbf{x}+\mathbf{t}^{-*}, \mathbf{y}-\mathbf{t}^{+*}\right) $ is
\begin{equation}
\label{eq:sbmbar}
\def\arraystretch{1}
\begin{array}{rl}
\rho ^*\left( \mathbf{x}+\mathbf{t}^{-*}, \mathbf{y}-\mathbf{t}^{+*}\right) =\min\limits_{\bm{\lambda},\mathbf{s}^-,\mathbf{s}^+} & \dfrac{1-\frac{1}{m}\sum _{i=1}^{m}s_i^- /\left( x_{i}+t_i^{-*}\right) }{1+\frac{1}{s}\sum _{r=1}^{s}s_r^+ /\left( y_{r}-t_r^{+*}\right)} \\
\text{s.t.} & \mathbf{x} +\mathbf{t}^{-*} = X\bm{\lambda}+\mathbf{s}^-, \\
& \mathbf{y} -\mathbf{t}^{+*} = Y\bm{\lambda}-\mathbf{s}^+, \\
& \bm{\lambda}\in \mathbb{R}^n_+,\quad \mathbf{s}^-\in \mathbb{R}^m_+,\quad \mathbf{s}^+\in \mathbb{R}^s_+.
\end{array}
\end{equation}

The following propositions allow us to simplify program \eqref{eq:sbmbar}. Moreover, Proposition \ref{propn2} implies that the programs given in \cite[Equation (2)]{Che19} and \cite[Equation (4)]{Che19} are equivalent.

\begin{proposition}
\label{propn}
Let $\mathbf{s}^{-*},\mathbf {s}^{+*}$ be optimal for \eqref{eq:sbmbar}, if $t_i^{-*}>0$ for some $i\in \left\{1,\ldots ,m\right\}$, then $s_i^{-*}=0$, and if $t_r^{+*}>0$ for some $r\in \left\{1,\ldots ,s\right\}$, then $s_r^{+*}=0$.
\end{proposition}

\begin{proposition}
\label{propn2} The objective function of \eqref{eq:sbmbar} can be replaced by
\begin{equation}
\label{eq:sbmbar2}
\frac{1-\frac{1}{m}\sum _{i=1}^{m}s_i^- /x_{i}}{1+\frac{1}{s}\sum _{r=1}^{s}s_r^+ /y_{r}}.
\end{equation}
\end{proposition}

According to \cite{Che19}, we define the \emph{continuous SBM} (CSBM) score (with respect to $\mathcal{D}$) of an activity $\left( \mathbf{x},\mathbf{y}\right) $ as
\begin{equation}
\label{eq:csbm}
\textrm{CSBM}\left( \mathbf{x},\mathbf{y}\right) =\dfrac{1-\frac{1}{m}\sum _{i=1}^{m}\left( s_i^{-*}-t_i^{-*}\right) /x_{i}}{1+\frac{1}{s}\sum _{r=1}^{s}\left( s_r^{+*}-t_r^{+*} \right) /y_{r}} ,
\end{equation}
where $\mathbf{t}^{-*},\mathbf{t}^{+*}$ are optimal super-efficiency slack vectors for \eqref{eq:ssbmact} and $\mathbf{s}^{-*},\mathbf{s}^{+*}$ are optimal inefficiency slack vectors for \eqref{eq:sbmbar}. Note that \eqref{eq:csbm} may not be well-defined if optimal slacks  $\mathbf{t}^{-*}$, $\mathbf{t}^{+*}$, $\mathbf{s}^{-*}$ or $\mathbf{s}^{+*}$ are not unique for the activity $\left( \mathbf{x},\mathbf{y}\right) $. In this case, the CSBM score may depend on which optimal slacks we choose.

The CSBM model can calculate both SBM efficiency (for activities in Region (I)) and SBM super-efficiency (for activities in Region (II)) scores, and it is indeed continuous. Nevertheless, as we show in Example \ref{chen19ex}, it may not be weakly monotonic for activities in Region (III) and hence, in these cases, it is not a valid score.

\begin{example}
\label{chen19ex}
We are going to consider a set of DMUs used in \cite{Doy93,Ton02} that consists of six efficient DMUs 
with four inputs and two outputs (see Table \ref{tab3data}).
In order to evaluate the first DMU, we are going to consider $\mathcal{D}=\left\{ \textrm{D2},\ldots ,\textrm{D6}\right\} $ as the set of reference DMUs. Then, the S-SBM score (with respect to $\mathcal{D}$) of D1 is $1.0116$ and its super-efficiency projection is $\left( \left( 80,627.89,54,8\right) ,\left( 90, 5\right) \right) $ whose SBM efficiency score is $0.7299$. The nonzero optimal slacks for programs \eqref{eq:ssbmact} and \eqref{eq:sbmbar} are $t_2^{-*}=27.89$, $s_4^{-*}=4.4$ and $s_2^{+*}=1.82$, giving a CSBM score of $0.7397$.

Now, let us consider $\textrm{D1}'$ equal to D1 except for the first input, that is increased in one unity, changing from $80$ to $81$. The S-SBM score (with respect to $\mathcal{D}$) of $\textrm{D1}'$ is $1.0103$ with super-efficiency projection $\left( \left( 81,624.82,54,8\right) ,\left( 90, 5\right) \right) $ whose SBM efficiency score is $0.7428$. The corresponding nonzero optimal slacks are $t_2^{-*}=24.82$, $s_4^{-*}=4.35$ and $s_2^{+*}=1.63$, giving a CSBM score of $0.7517$.

Since $\textrm{D1}'$ is strictly dominated by D1 and $\textrm{CSBM}\left( \textrm{D1}'\right) >\textrm{CSBM}\left( \textrm{D1}\right) $, we conclude that the CSBM score is not weakly monotonic in this case. Using the same technique as in the proofs of propositions \ref{prop:rhomon} and \ref{prop:22}, it can be proved that this type of example can only appear when an input or output that does not have any associated optimal slack is altered. In our case, when the first input is worsened, the decrease in the S-SBM score is not able to compensate for the increase in the SBM efficiency score of the super-efficiency projection.
\end{example}

\begin{table}[tbp]
\caption{Data of DMUs from Example \ref{chen19ex}.}
\begin{center}
\begin{tabular}{@{}lrrrrrr}
\toprule
& \multicolumn{4}{c}{Inputs} & \multicolumn{2}{c}{Outputs} \\
\cmidrule(lr){2-5} \cmidrule(lr){6-7}
\multicolumn{1}{c}{} & \multicolumn{1}{c}{$x_1$} & \multicolumn{1}{c}{$x_2$} & \multicolumn{1}{c}{$x_3$} & \multicolumn{1}{c}{$x_4$} & \multicolumn{1}{c}{$y_1$} & \multicolumn{1}{c}{$y_2$} \\
\midrule
D1 & 80 & 600 & 54 & 8 & 90 & 5 \\
D2 & 65 & 200 & 97 & 1 & 58 & 1 \\
D3 & 83 & 400 & 72 & 4 & 60 & 7 \\
D4 & 40 & 1000 & 75 & 7 & 80 & 10 \\
D5 & 52 & 600 & 20 & 3 & 72 & 8 \\
D6 & 94 & 700 & 36 & 5 & 96 & 6 \\
\bottomrule
\end{tabular}
\end{center}
\label{tab3data}
\end{table}

According to \cite{Fuk14}, there does not exist any weakly monotonic efficiency measure that uses a ratio-form least-distance approach to the closest projections over the efficient frontier. It seems that something similar can happen with the CSBM score, since \eqref{eq:csbm} is a ratio-form expression.

\section{The composite SBM model}
\label{sec:ceff0}

Example \ref{chen19ex} shows that optimal slacks (and hence the SBM efficiency score) of super-ef\-ficien\-cy projections do not serve to quantify the overestimation of the S-SBM score in some cases. Hence, for this purpose, we need  scores that do not just take into account the super-efficiency projection. Following this idea, in Subsection \ref{sec:ceff} we are going to define a continuous score function $\gamma $ that is equal to $\rho ^*$ in Region (I) and coincides with $\delta ^*$ in Region (II). Moreover, unlike the CSBM, $\gamma $ will always be weakly monotonic. Nevertheless, the computation of $\gamma $ involves nonlinear programming and hence, in Subsection \ref{sec:comp}, we are going to study some computational aspects.

\subsection{Definitions and properties}
\label{sec:ceff}

We define the \emph{composite SBM} (CompSBM) score (with respect to $\mathcal{D}$) of an activity $\left( \mathbf{x},\mathbf{y}\right) $ as
\begin{equation}
\label{eq:sbmcm}
\gamma \left( \mathbf{x},\mathbf{y}\right) =\delta ^*\left( \mathbf{x},\mathbf{y}\right) \cdot \max \rho ^*|_{\bar{P}} ,
\end{equation}
where $\bar{P}$ is the set of activities in $P$ that are dominated by $\left( \mathbf{x},\mathbf{y}\right) $ (see \eqref{eq:barP}) and $\max \rho ^*|_{\bar{P}}$ is the best (i.e. highest) SBM efficiency score of activities in $\bar{P}$, that can be interpreted as the SBM efficiency score of $\bar{P}$ as a set. Note that it is well-defined since $\bar{P}$ is closed.

The idea behind the CompSBM score $\gamma $ given by \eqref{eq:sbmcm} is not to focus only on super-efficiency projections, but on the entire set $\bar{P}$: instead of interpreting $\delta ^*\left( \mathbf{x},\mathbf{y}\right) $ as a distance from $\left( \mathbf{x},\mathbf{y}\right) $ to its super-efficiency projection and penalize the inefficiency of such projection, let us interpret $\delta ^*\left( \mathbf{x},\mathbf{y}\right) $ as a distance from $\left( \mathbf{x},\mathbf{y}\right) $ to $\bar{P}$ and penalize the inefficiency of $\bar{P}$, i.e. the fact that any activity in $\bar{P}$ is inefficient. These two points of view are not equivalent, as we will see in Remark \ref{rem:strmon}.

\begin{remark}[Unit-invariance]
\label{rem:comp}
In the J-SBM, CSBM and CompSBM models, we are implicitly assuming that slacks from different inputs and/or outputs can somehow compensate for each other, as pointed out and discussed by \cite{Che13}. In the case of the CompSBM model, according to expression \eqref{eq:sbmcm}, the optimal super-efficiency slacks of $\left( \mathbf{x},\mathbf{y}\right) $ (which are contained inside $\delta ^*\left( \mathbf{x},\mathbf{y}\right) $) are compensated by the optimal inefficiency slacks of the most efficient activity in $\bar{P}$ (which are contained inside $\max \rho ^*|_{\bar{P}}$). 
For this reason, unit-invariance is a very important property that these models should satisfy. On one hand, J-SBM and CSBM models are proved to be unit-invariant (see \cite{Che13,Che19}, respectively); on the other hand, the CompSBM model is unit-invariant because the SBM efficiency and super-efficiency models are unit-invariant (see \cite{Ton01,Ton02}).
\end{remark}

Next proposition clarifies the behavior of the CompSBM model, showing that it integrates the SBM efficiency $\rho ^*$ and the S-SBM $\delta ^*$ in its score.

\begin{proposition}
\label{prop41g2}
Let $\left( \mathbf{x},\mathbf{y}\right) $ be an activity and let $\bar{P}$ be the set given by \eqref{eq:barP}. Then
\begin{enumerate}
\item $\gamma \left( \mathbf{x},\mathbf{y}\right)  =\rho ^*\left( \mathbf{x},\mathbf{y}\right) $ if $\left( \mathbf{x},\mathbf{y}\right) $ is inefficient.
\item $\gamma \left( \mathbf{x},\mathbf{y}\right) =\delta ^*\left( \mathbf{x},\mathbf{y}\right) $ if $\left( \mathbf{x},\mathbf{y}\right) $ is efficient and $\max \rho ^*|_{\bar{P}}=1$, i.e. there are efficient activities in $\bar{P}$.
\item $\max \rho ^*|_{\bar{P}}<\gamma \left( \mathbf{x},\mathbf{y}\right) <\delta ^*\left( \mathbf{x},\mathbf{y}\right) $ if $\left( \mathbf{x},\mathbf{y}\right) $ is efficient and $\max \rho ^*|_{\bar{P}}<1$, i.e. any activity in $\bar{P}$ is inefficient.
\end{enumerate}
\end{proposition}

Analogously to what happens with the J-SBM and CSBM scores, Proposition \ref{prop41g2} ensures that 
$\gamma =\rho ^*$ in Region (I), $\gamma =\delta ^*$ in Region (II), and $\rho ^*\left( \bar{\mathbf{x}}^*,\bar{\mathbf{y}}^*\right) <\gamma \left( \mathbf{x},\mathbf{y}\right) \leq \delta ^*\left( \mathbf{x},\mathbf{y}\right) $ for $\left( \mathbf{x},\mathbf{y}\right) $ in Region (III), where $\left( \bar{\mathbf{x}}^*,\bar{\mathbf{y}}^*\right) $ is any super-efficiency projection of $\left( \mathbf{x},\mathbf{y}\right) $.
In the next results, we prove global continuity and weak monotonicity of the CompSBM score.

\begin{proposition}
\label{propcont}
The CompSBM score $\gamma $ is continuous.
\end{proposition}

\begin{proposition}
\label{prop42mon}
The CompSBM score $\gamma $ is weakly monotonic.
\end{proposition}

\begin{remark}[Super-inefficiency]
\label{rem2}
The CompSBM model is based on the SBM efficiency and super-efficiency models, providing a continuous score in the weakly efficient frontier.
However, since there is continuity, inevitably there will be efficient activities (in Region (III)) with CompSBM scores less than $1$ around the weakly efficient frontier, as it also happens with J-SBM and CSBM scores. Although this may seem a little counter-intuitive at first, \cite{Che13} gives a clarifying example in his ``Discussion and summary'' section.
Following the same criteria as \cite{Che13}, an efficient activity with score less than $1$ is said to be \emph{super-inefficient}.
In fact, \emph{super-inefficiency} is interpreted by \cite{Che13} as a ``hidden'' inefficiency, and \cite{Che19} affirms that super-inefficiency is a new division for efficiency, different from existing studies such as SBM efficiency and SBM super-efficiency. It is important to note that super-inefficiency has to be interpreted under the assumption that slacks from different inputs and/or outputs can somehow compensate for each other (see Remark \ref{rem:comp}). In the case of the CompSBM score, an efficient activity is super-inefficient if the model estimates that the magnitude of its optimal super-efficiency slacks is less than the magnitude of the optimal inefficiency slacks of the most efficient activity in $\bar{P}$.

Following \cite{Che13}, the \emph{super-inefficiency zone} is formed by all the super-inef\-ficient activities; on the other hand, the \emph{super-efficiency zone} is formed by all the efficient activities that are not in the super-inefficiency zone, i.e. with scores greater than or equal to 1. Each score (J-SBM, CSBM or CompSBM) can define a different super-inefficiency zone, but they are always in Region (III), around the weakly efficient frontier (see Figure \ref{figmapas}).
\end{remark}

According to Remark \ref{rem2}, any global continuous SBM model cannot determine whether an activity is efficient or not, because super-inefficient activities have scores less than $1$ but they are efficient. For this reason, any score that produces super-inefficient activities has to be taken as a complement to SBM efficiency and super-efficiency models. On the other hand, it could be interesting to define an alternative SBM super-efficiency score that penalizes the lack of efficient activities in $\bar{P}$ but it does not produce super-inefficient activities. Hence, in Remark \ref{rem:compssbm} we construct a super-efficiency score $\gamma_{\text{se}}$ based on $\gamma$ such that efficient activities obtain scores greater than or equal to $1$ (see \eqref{eq:gammase}). However, as it happens with the S-SBM score $\delta ^*$, discontinuities inevitably appear when implementing $\gamma_{\text{se}}$ in conjunction with SBM efficiency, leading to serious interpretation problems related to sensitivity.

\begin{remark}[Composite SBM super-efficiency score]
\label{rem:compssbm}
\cite{Lee21} defined a \emph{composite super-efficiency index} $\ddot \sigma  ^*$ such that efficient activities have scores greater than or equal to $1$ and the inefficiency of super-efficiency projections is penalized. Given an activity $ \left( \mathbf{x},\mathbf{y}\right) $, the corresponding score function (with respect to $\mathcal{D}$) would have this form:
\begin{equation}
\label{eq:ddotsigma}
\ddot \sigma  ^*\left( \mathbf{x},\mathbf{y}\right) =\left( \delta ^* \left( \mathbf{x},\mathbf{y}\right) -1\right) \cdot \rho ^* \left( \bar{\mathbf{x}}^*,\bar{\mathbf{y}}^*\right) +1,
\end{equation}
where $\left( \bar{\mathbf{x}}^*,\bar{\mathbf{y}}^*\right) $ is a super-efficiency projection of $ \left( \mathbf{x},\mathbf{y}\right) $. However, $\ddot \sigma ^*$ may not be well-defined if super-efficiency projections are not unique and, more importantly, it is not weakly monotonic in some cases, as we show in this example: considering the set of DMUs of Table \ref{tab3data}, we have that the S-SBM score of D2 (with respect to $\left\{ D1,D3,\ldots ,D6\right\} $) is $1.4146$ and the SBM efficiency of its super-efficiency projection is $0.3185$, giving a $\ddot \sigma ^*$ score of $1.132$; but increasing the fourth input of D2 from $1$ to $1.5$, we obtain that the S-SBM score is $1.3528$ and the SBM efficiency of the super-efficiency projection is $0.4392$, that gives a $\ddot \sigma ^*$ score of $1.1549$.

In order to fix this, based on \eqref{eq:sbmcm} and \eqref{eq:ddotsigma}, we define the \emph{composite SBM super-efficiency} (CompS-SBM) score (with respect to $\mathcal{D}$) of an activity $ \left( \mathbf{x},\mathbf{y}\right) $ by
\begin{equation}
\label{eq:gammase}
\begin{array}{rcl}
\gamma_{\text{se}}\left( \mathbf{x},\mathbf{y}\right) & = & \left( \delta ^* \left( \mathbf{x},\mathbf{y}\right) -1\right) \cdot \max \rho ^*|_{\bar{P}} +1 \\
& = & \gamma \left( \mathbf{x},\mathbf{y}\right) \cdot \left(1-1/\delta ^*\left( \mathbf{x},\mathbf{y}\right) \right) +1,
\end{array}
\end{equation}
which is always well-defined, unit-invariant, continuous and weakly monotonic. Moreover, it fulfils $\gamma_{\text{se}}=1$ in Region (I), $\gamma_{\text{se}}=\delta ^*$ in Region (II), and $1\leq \gamma_{\text{se}}\leq \delta ^*$ in Region (III), penalizing the lack of efficient activities in $\bar{P}$. It is important to remark that, although $\gamma_{\text{se}}$ fixes the overestimation of the S-SBM score, it obviously produces discontinuities in the weakly efficient frontier when implementing in conjunction with SBM efficiency, i.e. considering $\rho ^*$ for inefficient activities (Region (I)) and $\gamma_{\text{se}}$ for efficient ones (Region (II) and (III)).
\end{remark}

\begin{remark}[Strong monotonicity]
\label{rem:strmon}
We will show in Example \ref{ex3} (specifically with D5) that there exist uncommon cases in which super-efficiency projections are inefficient but there are efficient activities in $\bar{P}$. In these cases, $\max \rho ^*|_{\bar{P}}=1$ and hence, the CompSBM score $\gamma $ does not penalize the inefficiency of super-efficiency projections, unlike J-SBM and CSBM scores. As a consequence, $\gamma $ is not strongly monotonic for efficient activities, as it happens with the S-SBM score $\delta ^*$ (see Remark \ref{rem:ssbmnm}), although for other reasons and with a much lower frequency of cases.
\end{remark}

\begin{remark}[Alternative composite scores]
\label{rem:nomon}
We can use any weakly monotonic efficiency score $f$ instead of $\rho ^*$ in the definition of the CompSBM score $\gamma $ \eqref{eq:sbmcm}. In this case, we need continuity of $f|_P$ and the fulfilment of properties 1 and 2 of \cite{Fuk14} (see Section \ref{sec:2.2}). Hence, it is easy to prove that $\gamma $ results continuous and weakly monotonic, even if the efficiency score $f$ is not weakly monotonic. This makes it possible to use efficiency scores like the SBM-Max efficiency \citep{Ton10, Ton16}, that is not weakly monotonic in some cases \citep{Fuk14, And17}. It could be interesting since there is a close connection between SBM-Max efficiency and SBM super-efficiency models \citep{Ton17}.
\end{remark}

\subsection{Computational aspects}
\label{sec:comp}

Given an activity $\left( \mathbf{x},\mathbf{y}\right) $, the set $\bar{P}$ is formed by all the activities in $P$ that are dominated by $\left( \mathbf{x},\mathbf{y}\right) $ (see \eqref{eq:barP}). Hence,
\begin{equation}
\label{eq:rhobarP2}
\def\arraystretch{1}
\begin{array}{rl}
\max \rho ^*|_{\bar{P}}=\max \limits_{\bm{\lambda},\mathbf{t}^-,\mathbf{t}^+} &
\rho ^*\left( \mathbf{x}+\mathbf{t}^-,\mathbf{y}-\mathbf{t}^+\right) \\
\text{s.t.} & \mathbf{x}+\mathbf{t}^-\geq X\bm{\lambda}, \\
& \mathbf{0}< \mathbf{y}-\mathbf{t}^+\leq Y\bm{\lambda}, \\
& \bm{\lambda}\in \mathbb{R}^n_+,\quad \mathbf{t}^-\in \mathbb{R}^m_+,\quad \mathbf{t}^+\in \mathbb{R}^s_+,
\end{array}
\end{equation}
where, according to \eqref{eq:sbmact}, the objective function of \eqref{eq:rhobarP2} is
\begin{equation}
\label{eq:rhobar}
\def\arraystretch{1}
\begin{array}{rl}
\rho ^*\left( \mathbf{x}+\mathbf{t}^-,\mathbf{y}-\mathbf{t}^+\right) =\min\limits_{\bm{\lambda},\mathbf{s}^-,\mathbf{s}^+} & \rho  \left( \mathbf{x}+\mathbf{t}^-,\mathbf{y}-\mathbf{t}^+,\mathbf{s}^-,\mathbf{s}^+\right) \\
\text{s.t.} & \mathbf{x}+\mathbf{t}^-= X\bm{\lambda}+\mathbf{s}^-, \\
& \mathbf{y}-\mathbf{t}^+= Y\bm{\lambda}-\mathbf{s}^+, \\
& \bm{\lambda}\in \mathbb{R}^n_+,\quad \mathbf{s}^-\in \mathbb{R}^m_+,\quad \mathbf{s}^+\in \mathbb{R}^s_+.
\end{array}
\end{equation}
Note that $\bm{\lambda}$ in \eqref{eq:rhobarP2} and \eqref{eq:rhobar} are different internal variables of these programs. In fact, the constraints of \eqref{eq:rhobar} assure that the constraints of \eqref{eq:rhobarP2} involving $\bm{\lambda}$ are satisfied and hence, it suffices to demand nonnegativity of $\mathbf{t}^-,\mathbf{t}^+$ and $\mathbf{t}^+<\mathbf{y}$ in \eqref{eq:rhobarP2}. But, in this case, we have to note that \eqref{eq:rhobar} may result infeasible for some nonnegative small values of  $\mathbf{t}^-,\mathbf{t}^+$.
The next result allows us to simplify \eqref{eq:rhobarP2}.

\begin{proposition}
\label{prop71} Let $\left( \mathbf{x},\mathbf{y}\right) $ be an activity and let $\bar{P}$ be the set given by \eqref{eq:barP}. Then
\begin{equation}
\label{eq:maxmin}
\def\arraystretch{1}
\begin{array}{rrl}
\max \rho ^*|_{\bar{P}}= \max \limits_{\mathbf{t}^-,\mathbf{t}^+} & \min \limits_{\bm{\lambda},\mathbf{s}^-,\mathbf{s}^+} & \rho  \left( \mathbf{x},\mathbf{y},\mathbf{s}^-,\mathbf{s}^+\right) \\
& \text{s.t.} & \mathbf{x}+\mathbf{t}^- = X\bm{\lambda}+\mathbf{s}^-, \\
& & \mathbf{y}-\mathbf{t}^+ = Y\bm{\lambda}-\mathbf{s}^+, \\
& & \bm{\lambda}\in \mathbb{R}^n_+,\quad \mathbf{s}^-\in \mathbb{R}^m_+,\quad\mathbf{s}^+\in \mathbb{R}^s_+, \\
\text{s.t.} & \multicolumn{2}{l}{\mathbf{t}^-\in \mathbb{R}^m_+,\quad\mathbf{t}^+\in \mathbb{R}^s_+,\quad \mathbf{t}^+<\mathbf{y}.} 
\end{array}
\end{equation}
\end{proposition}

The inner minimization program of \eqref{eq:maxmin} can be linearized using the Charnes-Cooper transformation \citep{Cha62,Cha78}. Note that in the outer maximization program we only demand nonnegativity of $\mathbf{t}^-,\mathbf{t}^+$  and hence, the inner minimization program may result infeasible for some small values of  $\mathbf{t}^-,\mathbf{t}^+$. If we do not want this to happen, we must demand all the constraints of \eqref{eq:rhobarP2} in the outer program.
Note that program \eqref{eq:maxmin} can be viewed as a nonlinear maximization program, but it is also a continuous maximin problem with coupled constraints. Some methods for solving this kind of problems are provided by \cite{Shi80,Rus08,Tso09}, among others.

\begin{remark}[Lower bound]
\label{rem:low}
In some cases, computation of $\max \rho ^*|_{\bar{P}}$ by means of \eqref{eq:maxmin} may result too expensive (see Example \ref{ex3}). In these cases, 
we can compute a lower bound given by $\rho ^*\left( \bar{\mathbf{x}}^*,\bar{\mathbf{y}}^*\right) $ where $\left( \bar{\mathbf{x}}^*,\bar{\mathbf{y}}^*\right) $ is a super-efficiency projection of $\left( \mathbf{x},\mathbf{y}\right) $. Then,
\begin{equation}
\label{eq:gammalow}
\gamma _{\textrm{low}}\left( \mathbf{x},\mathbf{y}\right) =\delta ^*\left( \mathbf{x},\mathbf{y}\right)  \cdot \rho ^*\left( \bar{\mathbf{x}}^*,\bar{\mathbf{y}}^*\right) ,
\end{equation}
is a lower bound of the CompSBM score $\gamma \left( \mathbf{x},\mathbf{y}\right) $.
Usually, $\rho ^*\left( \bar{\mathbf{x}}^*,\bar{\mathbf{y}}^*\right) $ is very close (or even equal) to $\max \rho ^*|_{\bar{P}}$, and its computation only involves linear programs (applying the Charnes-Cooper transformation): \eqref{eq:ssbmact} for the super-efficiency projection and \eqref{eq:sbmbar} for its efficiency. On the other hand, \cite{Tso09} propose an algorithm for solving continuous maximin problems that requires a lower bound and hence, $\rho ^*\left( \bar{\mathbf{x}}^*,\bar{\mathbf{y}}^*\right) $ could serve for computing $\max \rho ^*|_{\bar{P}}$ using this algorithm.
\end{remark}

\section{Extensions}
\label{sec:ext}

In this section we are going to extend the CompSBM model to different orientations and returns to scale. Moreover, we discuss nonpositive data, and weighted inputs and/or outputs. Finally, we present a version adapted to the additive model.

\textbf{Orientations.} The CompSBM score is defined in a nonoriented form. Nevertheless, considering the input and output oriented versions of the S-SBM model \citep{Ton02}, we can adapt our CompSBM score to these orientations.
In this way, the \emph{input oriented CompSBM} score (with respect to $\mathcal{D}$) of an activity $\left( \mathbf{x},\mathbf{y}\right) $  is given by
$
\gamma _{\mathrm{I}}\left( \mathbf{x},\mathbf{y}\right) =\delta _{\mathrm{I}}^*\left( \mathbf{x},\mathbf{y}\right) \cdot \max \rho ^*_{\mathrm{I}}|_{\bar{P}}
$,
where $\delta ^*_{\mathrm{I}}$ and $\rho ^*_{\mathrm{I}}$ are the input oriented versions of $\delta ^*$ and $\rho ^*$ respectively. In order to compute $\max \rho ^*_{\mathrm{I}}|_{\bar{P}}$, it is easy to prove that
\begin{equation*}
\def\arraystretch{1}
\begin{array}{rrl}
1-\max \rho _{\textrm{I}}^*|_{\bar{P}}= \min \limits_{\mathbf{t}^-} & \max \limits_{\bm{\lambda},\mathbf{s}^-} & \frac{1}{m}\sum _{i=1}^{m}s_i^- /x_{i}  \\
& \text{s.t.} & \mathbf{x}+\mathbf{t}^- = X\bm{\lambda}+\mathbf{s}^-, \\
& & \bm{\lambda}\in \mathbb{R}^n_+,\quad \mathbf{s}^-\in \mathbb{R}^m_+, \\
\text{s.t.} & \multicolumn{2}{l}{\mathbf{t}^-\in \mathbb{R}^m_+,}
\end{array}
\end{equation*}
which is a linear continuous minimax problem with coupled constraints.
Analogously, the \emph{output oriented CompSBM} score (with respect to $\mathcal{D}$) of $\left( \mathbf{x},\mathbf{y}\right) $ is given by
$
\gamma _{\mathrm{O}}\left( \mathbf{x},\mathbf{y}\right) =\delta _{\mathrm{O}}^*\left( \mathbf{x},\mathbf{y}\right) \cdot \max \rho ^*_{\mathrm{O}}|_{\bar{P}}
$,
where $\delta ^*_{\mathrm{O}}$ and $\rho ^*_{\mathrm{O}}$ are the output oriented versions of $\delta ^*$ and $\rho ^*$ respectively. In order to compute $\max \rho ^*_{\mathrm{O}}|_{\bar{P}}$, it is easy to prove that
\begin{equation*}
\def\arraystretch{1}
\begin{array}{rrl}
1/\max \rho _{\textrm{O}}^*|_{\bar{P}}-1= \min \limits_{\mathbf{t}^+} & \max \limits_{\bm{\lambda},\mathbf{s}^+} & \frac{1}{m}\sum _{r=1}^{s}s_r^+ /y_{r}  \\
& \text{s.t.} & \mathbf{y}-\mathbf{t}^+ = Y\bm{\lambda}-\mathbf{s}^+, \\
& & \bm{\lambda}\in \mathbb{R}^n_+,\quad \mathbf{s}^+\in \mathbb{R}^s_+, \\
\text{s.t.} & \multicolumn{2}{l}{\mathbf{t}^+\in \mathbb{R}^s_+,\quad \mathbf{t}^+<\mathbf{y},}
\end{array}
\end{equation*}
which is also a linear continuous minimax problem with coupled constraints.

\textbf{Returns to scale.} In this work, all definitions and results are within the framework of constant returns to scale. Nevertheless, we can modify the programs to consider different returns to scale. For example, for variable returns to scale in the CompSBM model we have to add the constraint $\sum _{j=1}^n\lambda _j=1$ to programs \eqref{eq:ssbmact} and \eqref{eq:maxmin}.
It is worth noting that, taking variable returns to scale, nonoriented S-SBM models are always feasible \citep{Ton02}, but oriented models may result infeasible. In these cases, it will be impossible to compute oriented CompSBM scores.

\textbf{Zero or negative data.} The CompSBM model can accept zero or negative data as long as the SBM efficiency and super-efficiency models accept it. In fact, how to deal with zeros in data is discussed in \cite{Ton01, Ton02} and, more recently, how to handle with nonpositive data in general is discussed in \cite{Ton20,Lee21}.

\textbf{Weights.} We can consider different weights for each input and/or output. For example, we can compute the S-SBM score in \eqref{eq:ssbmact} by means of the weighted objective function
\begin{equation*}
\delta _w\left( \mathbf{x},\mathbf{y},\mathbf{t}^-,\mathbf{t}^+\right) = \dfrac{1+\frac{1}{\sum _{i=1}^m w_i^-}\sum _{i=1}^m w_i^- t_i^- /x_i}{1-\frac{1}{\sum _{r=1}^s w_r^+}\sum _{r=1}^s w_r^+ t_r^+ /y_r},
\end{equation*}
where $\mathbf{w}^-,\mathbf{w}^+$ are the corresponding weights vectors. Analogously, the SBM efficiency model would also have to take into account these weights.

\textbf{Additive model.} Finally, we can adapt the CompSBM score to the additive model. Following \cite{Cha82}, the \emph{additive efficiency} score (with respect to $\mathcal{D}$) of an activity $\left( \mathbf{x},\mathbf{y}\right) $ in unit-invariant form is defined as
\begin{equation}
\label{eq:addact}
\def\arraystretch{1}
\begin{array}{rl}
\alpha ^*\left( \mathbf{x},\mathbf{y}\right) =\max \limits_{\bm{\lambda},\mathbf{s}^-,\mathbf{s}^+} & \alpha \left( \mathbf{x},\mathbf{y},\mathbf{s}^-,\mathbf{s}^+\right) = \dfrac{1}{m+s}\left( \displaystyle \sum _{i=1}^m \dfrac{s_i^-}{x_i} +\sum _{r=1}^s \dfrac{s_r^+}{y_r}\right) \\
\text{s.t.} & \mathbf{x}= X\bm{\lambda}+\mathbf{s}^-, \\
& \mathbf{y}= Y\bm{\lambda}-\mathbf{s}^+, \\
& \bm{\lambda}\in \mathbb{R}^n_+,\quad \mathbf{s}^-\in \mathbb{R}^m_+,\quad\mathbf{s}^+\in \mathbb{R}^s_+.
\end{array}
\end{equation}
Note that $\alpha ^*$ is not an efficiency score satisfying properties $1$ and $2$ (see Section \ref{sec:2.2}). In fact, an activity $\left( \mathbf{x},\mathbf{y}\right) $ is efficient if and only if $\alpha ^*\left( \mathbf{x},\mathbf{y}\right) =0$.
On the other hand, following \cite{Du10}, we define the \emph{additive super-efficiency} score (with respect to $\mathcal{D}$) of $\left( \mathbf{x},\mathbf{y}\right) $ as
\begin{equation}
\label{eq:saddact}
\def\arraystretch{1}
\begin{array}{rl}
\beta ^*\left( \mathbf{x},\mathbf{y}\right) =\min \limits_{\bm{\lambda},\mathbf{t}^-,\mathbf{t}^+} & \beta \left( \mathbf{x},\mathbf{y},\mathbf{t}^-,\mathbf{t}^+\right) = \dfrac{1}{m+s}\left( \displaystyle \sum _{i=1}^m \dfrac{t_i^-}{x_i} +\sum _{r=1}^s \dfrac{t_r^+}{y_r}\right) \\
\text{s.t.} & \mathbf{x}+\mathbf{t}^- \geq X\bm{\lambda}, \\
& \mathbf{0}<\mathbf{y}-\mathbf{t}^+ \leq Y\bm{\lambda}, \\
& \bm{\lambda}\in \mathbb{R}^n_+,\quad \mathbf{t}^-\in \mathbb{R}^m_+,\quad\mathbf{t}^+\in \mathbb{R}^s_+.
\end{array}
\end{equation}
If $\left( \mathbf{x},\mathbf{y}\right) $ is inefficient, then $\beta ^*\left( \mathbf{x},\mathbf{y}\right) =0$.
So, from \eqref{eq:addact} and \eqref{eq:saddact} we can define the \emph{additive composite} score (with respect to $\mathcal{D}$) of $\left( \mathbf{x},\mathbf{y}\right) $ as
$
\gamma _{\textrm{add}} \left( \mathbf{x},\mathbf{y}\right) =\beta ^*\left( \mathbf{x},\mathbf{y}\right) -\min \alpha ^*|_{\bar{P}}
$.
Hence, $\gamma _{\textrm{add}}$ is negative for inefficient activities, and nonnegative for activities with efficient activities in $\bar{P}$.

\section{Examples}
\label{sec:examples}

In this section we are going to illustrate the J-SBM, CSBM and CompSBM models with some examples. 
We have used R 3.6.0 \citep{R19} for computations. Specifically, we have used the deaR package \citep{deaR19} for computing linear scores, and the NLopt package \citep{NLopt} for solving the nonlinear  program \eqref{eq:maxmin} in Example \ref{ex3}.

\begin{example}
\label{ex81}

Let us consider $\mathcal{D}=\left\{ \text{D1},\text{D2},\text{D3}\right\} $ a set of reference DMUs with two inputs and one normalized output, where $\text{D1}=\left( \left( 10,40\right) ,1\right) $, $\text{D2}=\left( \left( 15,25\right) ,1\right) $ and $\text{D3}=\left( \left( 30,20\right) ,1\right) $. In Figures \ref{figmapas} and \ref{figmapas_zoom}, we have computed binary logarithms of different scores of activities of the form $\left( \left( x_1,x_2\right) ,1\right) $. The super-inefficiency zones are represented, showing how discontinuity issues on weakly efficient activities are fixed by the CSBM and CompSBM scores, but not by the J-SBM score. Note that in cases where there is only one output (as in this example) the super-efficiency projections $\left( \bar{\mathbf{x}}^*,\bar{\mathbf{y}}^*\right) $ are among the most efficient activities in the corresponding set $\bar{P}$, and hence $\max \rho ^*|_{\bar{P}}=\rho ^*\left( \bar{\mathbf{x}}^*,\bar{\mathbf{y}}^*\right) $. Then, computing the CompSBM score $\gamma $ is equivalent to computing $\gamma _{\textrm{low}}$ (see \eqref{eq:gammalow}), whose program is linear.

\begin{figure}[tbp]
\begin{center}
\includegraphics[width=0.95\textwidth]{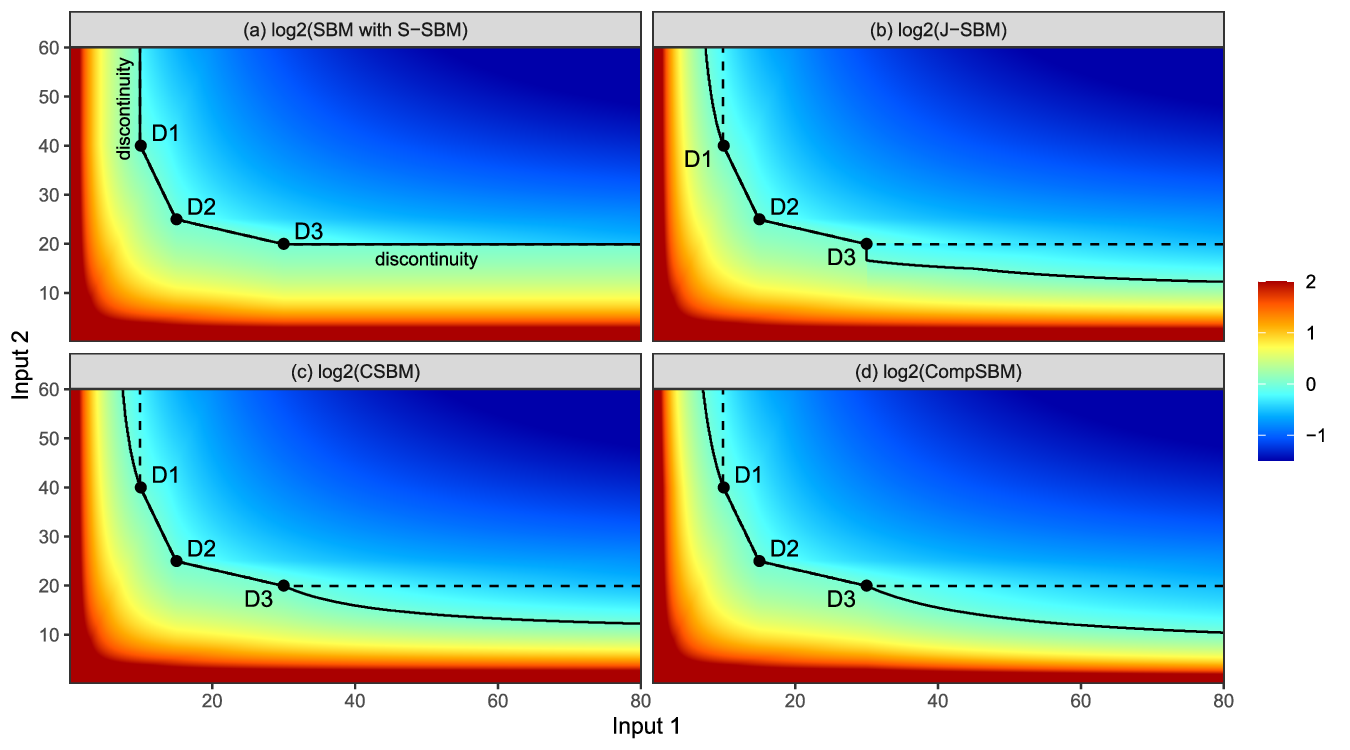}
\end{center}
\caption{Diagrams representing binary logarithms of different scores (with respect to the reference DMUs of Example \ref{ex81}) of activities with a normalized output equal to $1$: (a) SBM efficiency in conjunction with S-SBM, (b) J-SBM, (c) CSBM, (d) CompSBM. The solid lines separate activities with score $<1$ from activities with score $\geq 1$, and dashed lines represent the weakly efficient frontier. In plots (b), (c) and (d), the zones between solid and dashed lines are the super-inefficiency zones, fixing the discontinuity issues in (c) and (d).}
\label{figmapas}
\end{figure}

\begin{figure}[tbp]
\begin{center}
\includegraphics[width=0.95\textwidth]{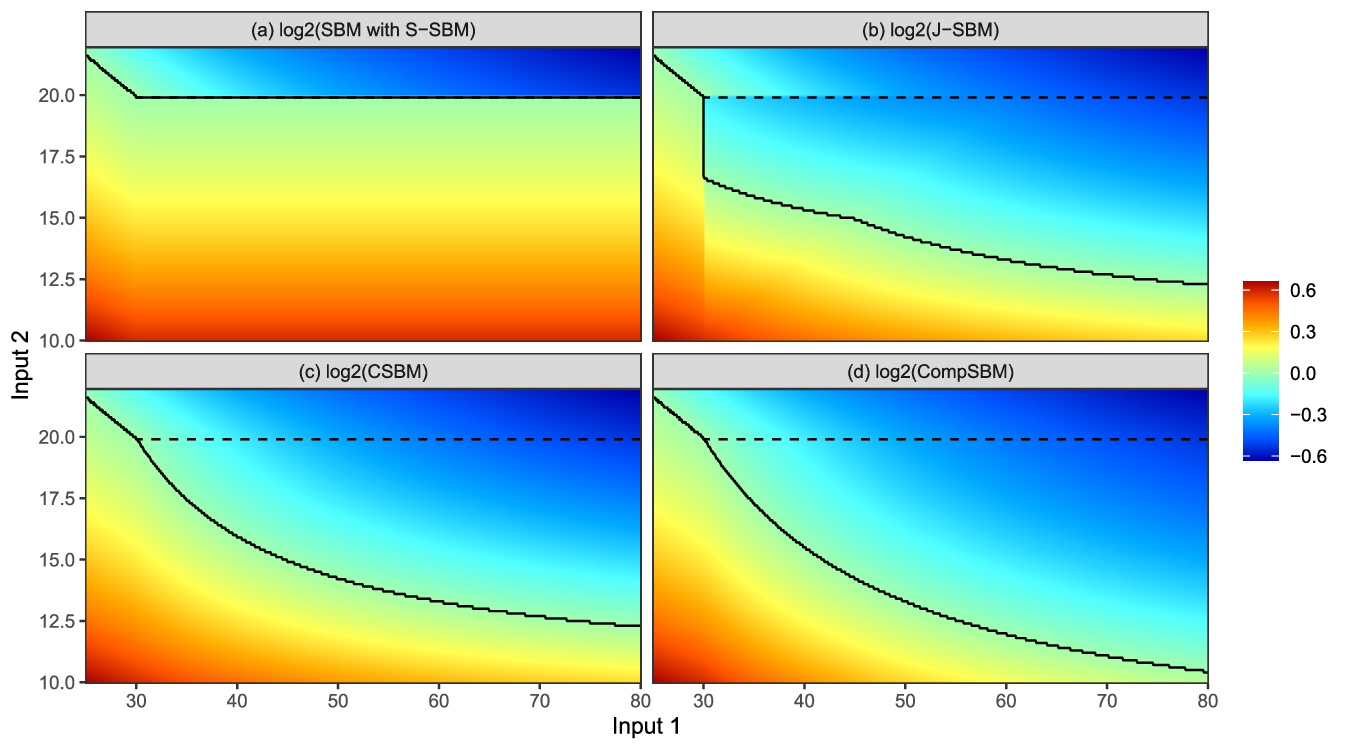}
\end{center}
\caption{Details of super-inefficiency zones in Figure \ref{figmapas}.}
\label{figmapas_zoom}
\end{figure}

Activities of the form $\left( \left( 30+c,x_2\right) ,1\right) $ with $0<x_2<20$ have the same S-SBM score $\delta ^*$ when $c\geq 0$ varies. But only when $c=0$ the super-efficiency projection is efficient and, in this case, it is equivalent to the fact that only when $c=0$ there are efficient activities in the corresponding $\bar{P}$. According to this, the J-SBM, CSBM and CompSBM models penalize activities of the form $\left( \left( 30+c,x_2\right) ,1\right) $ with $0<x_2<20$ and $c>0$.

\end{example}

\begin{example}
\label{exmi1}
In this example we want to illustrate how to introduce the CompSBM model into the Malmquist index computation. Note that this methodology is also applicable to the J-SBM and CSBM models. 
We consider the set of DMUs $\mathcal{D}=\left\{ \text{D1},\text{D2},\text{D3}\right\} $ of Example \ref{ex81} and another DMU D4 with activity changing from $\text{A1}=\left( \left( 30,10\right) ,1\right) $ to $\text{A2}=\left( \left( 40,10\right) ,1\right) $, while the activities of the DMUs in $\mathcal{D}$ remain unchanged. In this way, using the SBM efficiency and super-efficiency models, we can compute the \emph{SBM Malmquist index} \citep{Ton04} of D4, that is the product of two factors: the \emph{catch-up} and the \emph{frontier-shift}. On one hand, the catch-up (or recovery) is interpreted as the DMU's relative efficiency change, so catch-up values greater than $1$ indicate progress, values less than $1$ indicate regress, and a catch-up equal to $1$ means no change. On the other hand, the frontier-shift (or innovation) is related to the technological change in the efficient frontiers and hence, analogously to the catch-up, values greater, less and equal to $1$ indicate progress, regress, and no change, respectively, of the efficient frontier with respect to the evaluated DMU.

Nevertheless, we can also compute a Malmquist index using the CompSBM score (or the J-SBM or CSBM scores) instead of the S-SBM score. Table \ref{tabex81-2} shows the scores of A1 and A2, while Table \ref{tab1c} shows the results of the Malmquist index with no orientation and using the exclusive scheme (see \cite{Ton04} for details). The original SBM Malmquist index does not detect the catch-up, and its frontier-shift does not take into account the inefficiency of the super-efficiency projection of A2. On the other hand, these problems do not appear when the CompSBM score is used instead, obtaining a catch-up smaller than one and a frontier-shift indicating less technological regress of the efficient frontier with respect to D4.

\begin{table}[tbp]
\caption{SBM efficiency score ($\rho ^*$), S-SBM score ($\delta ^*$), SBM efficiency score of super-efficiency projections ($\rho ^*\left( \bar{\mathbf{x}}^*,\bar{\mathbf{y}}^*\right) $) that coincides with the best SBM efficiency score in the corresponding $\bar{P}$ ($\max \rho ^*|_{\bar{P}}$), and CompSBM score ($\gamma $) of activities $\text{A1},\ldots ,\text{A4}$ from Example \ref{exmi1}.}
\begin{center}
\begin{tabular}{@{}rrrrrr}
\toprule
& \multicolumn{1}{c}{$\rho ^*$} & \multicolumn{1}{c}{$\delta ^*$} & \multicolumn{1}{c}{$\rho ^*\left( \bar{\mathbf{x}}^*,\bar{\mathbf{y}}^*\right) $} & \multicolumn{1}{c}{$\max \rho ^*|_{\bar{P}}$} & \multicolumn{1}{c}{$\gamma $} \\
\midrule
A1  & 1 & 1.5 & 1& 1 & 1.5 \\
A2  & 1 & 1.5 & 0.875 & 0.875 & 1.3125 \\
A3  & 0.8 & 1 & 0.8 & 0.8 & 0.8 \\
A4  & 1 & 1.0263 & 0.8 & 0.8 & 0.821 \\
\bottomrule
\end{tabular}
\end{center}
\label{tabex81-2}
\end{table}

\begin{table}[tbp]
\caption{Catch-up, frontier-shift and Malmquist index of DMU D4 changing from A1 to A2 and from A3 to A4 (see Example \ref{exmi1}) in the original SBM form and using the CompSBM score, with no orientation and taking the exclusive scheme.}
\begin{center}
\begin{tabular}{@{}clrrr}
\toprule
& & \multicolumn{1}{c}{Catch-up} & \multicolumn{1}{c}{Frontier-shift} & \multicolumn{1}{c}{Malmquist index} \\
\midrule
\multirow{2}{*}{$\text{A1}\rightarrow \text{A2}$} & SBM & 1 & 0.866 &  0.866 \\
& CompSBM & 0.875 & 0.9258 & 0.8101 \\
\midrule
\multirow{2}{*}{$\text{A3}\rightarrow \text{A4}$} & SBM  & 1.2829 & 1 &  1.2829 \\
& CompSBM  & 1.0262 & 1 & 1.0262 \\
\bottomrule
\end{tabular}
\end{center}
\label{tab1c}
\end{table}

Now, let us consider that the activity of D4 changes from $\text{A3}=\left( \left( 50,20\right) ,1\right) $ (that is weakly efficient) to $\text{A4}=\left( \left( 50,19\right) ,1\right) $, while the activities of the DMUs in $\mathcal{D}$ remain unchanged. Analogously, we can compute the SBM Malmquist index of D4 in its original form or using the CompSBM score instead. Table \ref{tabex81-2} shows the scores of A3 and A4, while Table \ref{tab1c} also shows the results of the Malmquist index with no orientation and using the exclusive scheme. The catch-up of the original SBM Malmquist index is overestimated due to the discontinuity of the SBM efficiency score on weakly efficient activities described in Example \ref{ex81}. This problem is again solved when the CompSBM score is used instead, resulting in a catch-up very close to $1$. 

\end{example}

\begin{example}
\label{ex3}

We are going to consider the same set of DMUs of Example \ref{chen19ex}, used in \cite{Doy93,Ton02}, consisting of six efficient DMUs (power plant locations) with four inputs and two outputs (see Table \ref{tab3data}).
Table \ref{tab3no} shows nonoriented scores (although \cite{Ton02} only considers the input oriented scenario) and Table \ref{tab3ets} shows the corresponding efficient targets of super-efficiency projections, with optimal slacks in parentheses. Note that the difference between the original DMU and the efficient target of a super-efficiency projection is given either by optimal super-efficiency slacks or by optimal inefficiency slacks (of the super-efficiency projection), separately in each input and output. This fact is shown in Table \ref{tab3ets}, where super-efficiency slacks are displayed in green, inefficiency slacks in red, and each input or output has either green or red slacks, but not both (see Proposition \ref{propn}). It should be also noted that the S-SBM score only takes into account ``green slacks'', ignoring ``red slacks''.

\begin{table}[tbp]
\caption{Nonoriented scores of DMUs from Example \ref{ex3}. The scores are as follows: S-SBM  score ($\delta ^*$), SBM efficiency score of super-efficiency projections ($\rho ^*\left( \bar{\mathbf{x}}^*,\bar{\mathbf{y}}^*\right) $), the best SBM efficiency score in the corresponding $\bar{P}$ ($\max \rho ^*|_{\bar{P}}$), J-SBM score, CSBM score, a lower bound of the CompSBM score ($\gamma _{\text{low}}=\delta ^* \cdot \rho ^*\left( \bar{\mathbf{x}}^*,\bar{\mathbf{y}}^*\right) $, see Remark \ref{rem:low}), CompSBM score ($\gamma $), and CompS-SBM score ($\gamma _{\textrm{se}}$), all of them computed with respect to the set of reference DMUs given by all the DMUs excluding the evaluated DMU.}
\begin{center}
\begin{tabular}{@{}lrrrrrrrr}
\toprule
DMU & \multicolumn{6}{l}{Scores} \\ 
\cmidrule(lr){2-9}
& \multicolumn{1}{c}{$\delta ^*$} & \multicolumn{1}{c}{$\rho ^*\left( \bar{\mathbf{x}}^*,\bar{\mathbf{y}}^*\right) $} & \multicolumn{1}{c}{$\max \rho ^*|_{\bar{P}}$} & \multicolumn{1}{c}{J-SBM} & \multicolumn{1}{c}{CSBM} & \multicolumn{1}{c}{$\gamma _{\textrm{low }}$} & \multicolumn{1}{c}{$\gamma $} & \multicolumn{1}{c}{$\gamma _{\textrm{se}}$} \\
\midrule
D1 & 1.0116 & 0.7299 & 0.8565 & 0.4990 & 0.7397 & 0.7384 & 0.8664 & 1.0099 \\
D2 & 1.4146 & 0.3185 & 0.4525 & 0.3791 & 0.3791 & 0.4505 & 0.6401 & 1.1876 \\
D3 & 1.0781 & 0.597 & 0.597 & 0.6732 & 0.6732 & 0.6436 & 0.6436 & 1.0466 \\
D4 & 1.1563 & 0.6162 & 0.6162 & 0.7633 & 0.7633 & 0.7125 & 0.7125 & 1.0963 \\
D5 & 1.5859 & 0.9687 & 1 & 1.5458 & 1.5458 & 1.5363 & 1.5859 & 1.5859 \\
D6 & 1.0198 & 0.7422 & 0.8227 & 0.6158 & 0.7577 & 0.7569 & 0.839 &  1.0163 \\
\bottomrule
\end{tabular}
\end{center}
\label{tab3no}
\end{table}

\begin{table}[tbp]
\caption{Efficient targets of super-efficiency projections of DMUs from Example \ref{ex3} with optimal slacks in parentheses. Super-efficiency slacks $t^{-*}_i,t^{+*}_r$ are displayed in green and inefficiency slacks (of super-efficiency projections) $s^{-*}_i,s^{+*}_r$ in red, for $i=1,\ldots ,4$ and $r=1,2$.}
\begin{center}
\begin{tabular}{@{}lllllll}
\toprule
DMU & \multicolumn{6}{l}{Efficient targets (and slacks)} \\ 
\cmidrule(lr){2-7}
& \multicolumn{1}{c}{$x_1$} & \multicolumn{1}{c}{$x_2$} & \multicolumn{1}{c}{$x_3$} & \multicolumn{1}{c}{$x_4$} & \multicolumn{1}{c}{$y_1$} & \multicolumn{1}{c}{$y_2$} \\
\midrule
D1 & 80 & 627.89 ({\color{OliveGreen} 27.89}) & 54 & 3.6 ({\color{red} 4.4}) & 90 & 6.82 ({\color{red} 1.82}) \\
D2 & 17.3 ({\color{red} 47.7}) & 200 & 6.7 ({\color{red} 90.3}) & 1 & 24 ({\color{OliveGreen} 34}) & 2.67 ({\color{red} 1.67}) \\
D3 & 45.5 ({\color{red} 37.5}) & 525 ({\color{OliveGreen} 125}) & 17.5 ({\color{red} 54.5}) & 2.62 ({\color{red} 1.38}) & 63 ({\color{red} 3}) & 7 \\
D4 & 65 ({\color{OliveGreen} 25}) & 750 ({\color{red} 250}) & 25 ({\color{red} 50}) & 3.75 ({\color{red} 3.25}) & 90 ({\color{red} 10}) & 10 \\
D5 & 70.5 ({\color{OliveGreen} 18.5}) & 525 ({\color{red} 75}) & 27 ({\color{OliveGreen} 7}) & 3.75 ({\color{OliveGreen} 0.75}) & 72 & 4.5 ({\color{OliveGreen} 3.5}) \\
D6 & 73.6 ({\color{red} 20.4}) & 755.5 ({\color{OliveGreen} 55.5}) & 36 & 5 & 96 & 9.29 ({\color{red} 3.29}) \\
\bottomrule
\end{tabular}
\end{center}
\label{tab3ets}
\end{table}

Note that, for D5, $\rho ^*\left( \bar{\mathbf{x}}^*,\bar{\mathbf{y}}^*\right) < 1$ and it does not coincide with $\max \rho ^*|_{\bar{P}}$, whose value is $1$. Hence, D5 has super-efficiency projections that are inefficient, but there are also efficient activities in the corresponding $\bar{P}$. In this case, the CompSBM model does not penalize the inefficiency of such super-efficiency projections (because there are efficient activities in $\bar{P}$), contrary to the J-SBM and CSBM models. An important conclusion is that the CompSBM score is not strongly monotonic in some cases, because D5 can be improved (specifically, the second input can be lowered up to $525$) keeping the original S-SBM score and, since there are efficient activities in $\bar{P}$, the CompSBM score of D5 will not change.

We have used the NLopt package \citep{NLopt} for solving the nonlinear program \eqref{eq:maxmin} in the computation of $\max \rho ^*|_{\bar{P}}$. Specifically, we have used the following global nonlinear algorithms included in this package: DIRECT (Dividing RECtangles) \citep{Jon93}, its ``locally biased'' version DIRECT-L \citep{Gab01}, and COBYLA (Constrained Optimization BY Linear Approximations) \citep{Pow98}. Computation time varies depending on the algorithm: 150 seconds for DIRECT, 20 seconds for DIRECT-L, and 2 seconds for COBYLA, using a 2 GHz processor. Nevertheless, computation time grows exponentially with the number of efficient DMUs, inputs and/or outputs, making it practically impossible to solve problems with, for example, more than 30 efficient DMUs with 5 inputs and 5 outputs.

\end{example}

\section{Concluding remarks}
\label{sec:conc}

The problem of ranking efficient DMUs continues to be an active issue that keeps generating new studies and methods within DEA \citep{Jab12, Zyk22}.
Since its inception in the papers of \cite{Pas99} and \cite{Ton01}, the SBM super-efficiency model has proved to be an important tool for ranking efficient DMUs which has been widely used by DEA practitioners in the last years. However, some interpretation problems appear when super-efficiency projections are weakly efficient, since the SBM super-efficiency model does not take into account the inefficiency slacks of such projections and therefore, some efficient DMUs can not be properly ranked. Moreover, this fact is closely related to discontinuities produced in the weakly efficient frontier when implementing SBM super-efficiency in conjunction with SBM efficiency.
Authors like \cite{Che13} and \cite{Che19} tried to solve these problems, but they did not arrive to a fully satisfactory solution. Nevertheless, their papers lay the foundations for future studies on the subject.

In order to shed some light on this matter, we have introduced the CompSBM model, leading to the first example of a weakly monotonic score that integrates the SBM efficiency and S-SBM scores in a continuous way. Indeed, we have shown that it coincides with the S-SBM score when the DMUs are efficient with no inefficient super-efficiency projections, and coincides with the SBM efficiency score when the DMUs are inefficient. Moreover, the CompSBM model gives a continuous ranking of DMUs, avoiding the abrupt changes in the scores shown by other models and hence, solving the discontinuity problems in the weakly efficient frontier. It can also be adapted to other models such as the additive model or the SBM Malmquist index, or even we can use alternative efficiency scores different from $\rho ^*$ in the construction of a composite score, giving us a deeper insight into the evolution of the performance of DMUs.

The idea behind the CompSBM model is not to consider only the super-efficiency slacks, but also to penalize the lack of efficient activities in the set $\bar{P}$ onto which super-efficiency projections are projected. Since we demand continuity, super-inefficient activities (i.e. efficient activities with scores less than $1$) inevitably appear around the weakly efficient frontier. Although it can be counter-intuitive at first, super-inefficiency is interpreted as a ``hidden'' inefficiency \citep{Che13}, assuming that slacks from different inputs/outputs can somehow compensate for each other. According to \cite{Che19}, super-inefficiency is a new division for efficiency, different from existing studies such as SBM efficiency and SBM super-efficiency. Nevertheless, since super-inefficiency is a relatively new concept, it is not yet considered by some researchers who prefer to deal with discontinuities rather than with super-inefficiencies, although discontinuities have serious interpretation problems related to sensitivity. Hence, we have also defined a new weakly monotonous SBM super-efficiency score (based on the CompSBM score and the work of \cite{Lee21}) that penalizes the lack of efficient activities in $\bar{P}$ without producing super-inefficient activities (see Remark \ref{rem:compssbm}). However, discontinuities obviously appear in the weakly efficient frontier when implementing this new super-efficiency in conjunction with SBM efficiency.

To sum up, the CompSBM score:
\begin{itemize}
\item is continuous,
\item is weakly monotonous, and
\item it allows ranking the efficient DMUs.
\end{itemize}
However, this score presents two difficulties, namely: 
\begin{itemize}
\item its calculation requires solving nonlinear optimization problems, and
\item it presents super-inefficiencies, which is the price we have to pay for having a global continuous score. 
\end{itemize}

We believe that the methodology employed in the construction of the CompSBM score can help in the development of other models with better properties, as for example to be strongly monotonic (see Remark \ref{rem:strmon}) or be easier to compute. Moreover, there are some pending tasks that may be interesting, such as the use of the SBM-Max efficiency model in the construction of the composite score (see Remark \ref{rem:nomon}), or the study of the potential infeasibility of oriented composite models under variable returns to scale. In our opinion, these and other questions could lead to future results which, beyond any doubt, will help to increase the understanding of the DEA methodology.

\section*{Declarations}

\begin{itemize}
\item[] Funding: Not applicable.

\item[] Conflicts of interest/Competing interests: Not applicable.

\item[] Availability of data and material: Not applicable.

\item[] Code availability: We have used R 3.6.0 \citep{R19} for computations. Specifically, we have used the deaR package \citep{deaR19} for computing linear scores, and the NLopt package \citep{NLopt} for solving the nonlinear  program \eqref{eq:maxmin} in Example \ref{ex3}.
\end{itemize}

\appendix
\section{Appendix}
\label{sec:app}

\begin{proof}[\textbf{Proof of Proposition \ref{prop:rhomon}}]
We are going to prove that $\dfrac{\partial \rho ^*}{\partial x_k} <0$ for $k\in \left\{1,\ldots ,m\right\} $ and $\dfrac{\partial \rho ^*}{\partial y_k} >0$ for $k\in \left\{1,\ldots ,s\right\} $ using the Karush-Kuhn-Tucker (KKT) conditions.
Let $\left( \mathbf{x},\mathbf{y}\right) $ be an activity in $P$ and $\mathbf{s}^{-*},\mathbf{s}^{+*}$ be optimal inefficiency slack vectors for program \eqref{eq:sbmact}. Given $k\in \left\{1,\ldots ,m\right\} $, we replace in \eqref{eq:sbmact} the parameter $x_k$ by a new variable $x$ and add the constraint $x=x_k$ whose KKT multiplier is $\dfrac{\partial \rho ^*}{\partial x_k}\left( \mathbf{x},\mathbf{y}\right) $. Hence, the stationarity condition associated to the new variable $x$ is given by
\begin{equation}
\label{eq:p1-1}
\frac{\partial \rho ^*}{\partial x_k}\left( \mathbf{x},\mathbf{y}\right) = \frac{\partial \rho }{\partial x_k}\left( \mathbf{x},\mathbf{y},\mathbf{s}^{-*},\mathbf{s}^{+*}\right) -\mu ^-_k ,
\end{equation}
where $\mu ^-_k$ is the KKT multiplier of the $k$-th constraint in \eqref{eq:sbmact}. The stationarity condition associated to $s^-_k$ of \eqref{eq:sbmact} is given by
\begin{equation}
\label{eq:p1-2}
\frac{\partial \rho }{\partial s^-_k}\left( \mathbf{x},\mathbf{y},\mathbf{s}^{-*},\mathbf{s}^{+*}\right) +\mu ^-_k=\nu[s^-_k]\geq 0,
\end{equation}
where $\nu[s^-_k]$ is the KKT multiplier of the nonnegativity condition of $s^-_k$. Hence, applying \eqref{eq:p1-2} in \eqref{eq:p1-1} we obtain
\begin{equation*}
\label{eq:p1-3}
\begin{array}{rcl}
\dfrac{\partial \rho ^*}{\partial x_k} \left( \mathbf{x},\mathbf{y}\right) & \leq & \dfrac{\partial \rho }{\partial x_k}\left( \mathbf{x},\mathbf{y},\mathbf{s}^{-*},\mathbf{s}^{+*}\right)+\dfrac{\partial \rho }{\partial s^-_k}\left( \mathbf{x},\mathbf{y},\mathbf{s}^{-*},\mathbf{s}^{+*}\right) \\[10pt]
& = & \dfrac{1}{1+\frac{1}{s}\sum _{r=1}^{s}s_r^{+*} /y_{r}}\dfrac{1}{mx_k}\left( \dfrac{s^{-*}_k}{x_k}-1\right) <0,
\end{array}
\end{equation*}
since $s^{-*}_k<x_k$.
Analogously, given $k\in \left\{1,\ldots ,s\right\} $, it can be proved that
\begin{equation*}
\label{eq:p1-4}
\begin{array}{rcl}
\dfrac{\partial \rho ^*}{\partial y_k} \left( \mathbf{x},\mathbf{y}\right) & \geq & \dfrac{\partial \rho }{\partial y_k}\left( \mathbf{x},\mathbf{y},\mathbf{s}^{-*},\mathbf{s}^{+*}\right)-\dfrac{\partial \rho }{\partial s^+_k}\left( \mathbf{x},\mathbf{y},\mathbf{s}^{-*},\mathbf{s}^{+*}\right) \\[10pt]
& = & \dfrac{\rho \left( \mathbf{x},\mathbf{y},\mathbf{s}^{-*},\mathbf{s}^{+*}\right) }{1+\frac{1}{s}\sum _{r=1}^{s}s_r^{+*} /y_{r}}\dfrac{1}{sy_k}\left( \dfrac{s^{+*}_k}{y_k}+1\right) >0.
\end{array}
\end{equation*}
\end{proof}

\begin{proof}[\textbf{Proof of Proposition \ref{prop:22}}]
Let $\left( \mathbf{x},\mathbf{y}\right) $ be an activity and $\mathbf{t}^{-*},\mathbf{t}^{+*}$ be optimal super-ef\-fi\-cien\-cy slack vectors for program \eqref{eq:ssbmact}. Given $k\in \left\{1,\ldots ,m\right\} $, we replace in \eqref{eq:ssbmact} the parameter $x_k$ by a new variable $x$ and add the constraint $x=x_k$ whose KKT multiplier is $\dfrac{\partial \delta ^*}{\partial x_k}\left( \mathbf{x},\mathbf{y}\right) $. Hence, the stationarity condition associated to the new variable $x$ is given by
\begin{equation}
\label{eq:p2-1}
\frac{\partial \delta ^*}{\partial x_k}\left( \mathbf{x},\mathbf{y}\right) = \frac{\partial \delta }{\partial x_k}\left( \mathbf{x},\mathbf{y},\mathbf{t}^{-*},\mathbf{t}^{+*}\right) -\mu ^-_k ,
\end{equation}
where $\mu ^-_k$ is the KKT multiplier of the $k$-th constraint in \eqref{eq:ssbmact}.

Let us suppose that $t_k^{-*}>0$. The KKT stationarity condition associated to $t^-_k$ of \eqref{eq:ssbmact} is given by
\begin{equation}
\label{eq:p2-2}
\frac{\partial \delta }{\partial t^-_k}\left( \mathbf{x},\mathbf{y},\mathbf{t}^{-*},\mathbf{t}^{+*}\right) -\mu ^-_k=\nu[t^-_k]= 0,
\end{equation}
where $\nu[t^-_k]$ is the KKT multiplier of the nonnegativity condition of $t^-_k$, that vanishes in virtue of the corresponding complementary slackness condition. Hence, applying \eqref{eq:p2-2} in \eqref{eq:p2-1} we obtain
\begin{equation}
\label{eq:p2-3}
\begin{array}{rcl}
\dfrac{\partial \delta ^*}{\partial x_k} \left( \mathbf{x},\mathbf{y}\right) & = & \dfrac{\partial \delta }{\partial x_k}\left( \mathbf{x},\mathbf{y},\mathbf{t}^{-*},\mathbf{t}^{+*}\right)-\dfrac{\partial \delta }{\partial t^-_k}\left( \mathbf{x},\mathbf{y},\mathbf{t}^{-*},\mathbf{t}^{+*}\right) \\[10pt]
& = & \dfrac{-1}{1-\frac{1}{s}\sum _{r=1}^{s}t_r^{+*} /y_{r}}\dfrac{1}{mx_k}\left( \dfrac{t^{-*}_k}{x_k}+1\right) <0,
\end{array}
\end{equation}
since $t_r^{+*}<y_r$ for $r=1,\ldots ,s$.

On the other hand, let us suppose that $t_k^{-*}=0$. From the KKT dual feasibility conditions of \eqref{eq:ssbmact}, we have $\mu ^-_k\geq 0$ and hence, from \eqref{eq:p2-1} we obtain
\begin{equation}
\label{eq:p2-4}
\frac{\partial \delta ^*}{\partial x_k}\left( \mathbf{x},\mathbf{y}\right) \leq \frac{\partial \delta }{\partial x_k}\left( \mathbf{x},\mathbf{y},\mathbf{t}^{-*},\mathbf{t}^{+*}\right) =\dfrac{-1}{1-\frac{1}{s}\sum _{r=1}^{s}t_r^{+*} /y_{r}}\dfrac{t^{-*}_k}{mx_k^2} = 0.
\end{equation}

Analogously to \eqref{eq:p2-3} and \eqref{eq:p2-4}, it can be proved that given $k\in \left\{1,\ldots ,s\right\} $, if $t_k^{+*}>0$ then $\dfrac{\partial \delta ^*}{\partial y_k} \left( \mathbf{x},\mathbf{y}\right) >0$, and if $t_k^{+*}=0$ then $\dfrac{\partial \delta ^*}{\partial y_k} \left( \mathbf{x},\mathbf{y}\right) \geq 0$. 
\end{proof}

\begin{proof}[\textbf{Proof of Proposition \ref{propn}}]
We define $\bm{\tau}^-,\bm{\tau}^+$ such that $\tau_i^-=\max \left\{ 0,t_i^{-*}-s_i^{-*}\right\} $, $\tau_r^+=\max \left\{ 0,t_r^{+*}-s_r^{+*}\right\} $, for $i=1,\ldots ,m$ and $r=1,\ldots ,s$.
From the constraints of \eqref{eq:sbmbar}, we have
\begin{equation}
\label{eq:propn1}
\begin{array}{l}
\mathbf{x} +\bm{\tau}^- \geq \mathbf{x} +\left( \mathbf{t}^{-*}-\mathbf{s}^{-*}\right) = X\bm{\lambda}^*, \\
\mathbf{y} -\bm{\tau}^+ \leq\mathbf{y} -\left( \mathbf{t}^{+*}-\mathbf{s}^{+*}\right) = Y\bm{\lambda}^*,
\end{array}
\end{equation}
where $\bm{\lambda}^*$ is optimal for \eqref{eq:sbmbar}, along with $\mathbf{s}^{-*},\mathbf {s}^{+*}$.
From \eqref{eq:propn1} and taking into account that $\bm{\tau}^-,\bm{\tau}^+$ are nonnegative, we have that $\left( \bm{\lambda}^*, \bm{\tau}^-, \bm{\tau}^+\right) $ satisfies the constraints of \eqref{eq:ssbmact}, and hence
\begin{equation}
\label{eq:propn2}
\delta \left( \mathbf{x},\mathbf{y},\bm{\tau}^-,\bm{\tau}^+\right) \geq \delta \left( \mathbf{x},\mathbf{y},\mathbf{t}^{-*},\mathbf{t}^{+*}\right) .
\end{equation}
On the other hand, if we suppose that there exists $i\in \left\{1,\ldots ,m\right\}$ such that $t_i^{-*},s_i^{-*}>0$, then $\tau _i^-<t_i^{-*}$ and hence $\delta \left( \mathbf{x},\mathbf{y},\bm{\tau}^-,\bm{\tau}^+\right) < \delta \left( \mathbf{x},\mathbf{y},\mathbf{t}^{-*},\mathbf{t}^{+*}\right) $, that contradicts \eqref{eq:propn2}. Analogously, if we suppose that there exists $r\in \left\{1,\ldots ,s\right\}$ such that $t_r^{+*},s_r^{+*}>0$, then we arrive to the same contradiction.
\end{proof}

\begin{proof}[\textbf{Proof of Proposition \ref{propn2}}]
Let \eqref{eq:sbmbar}$'$ be the program \eqref{eq:sbmbar} with objective function \eqref{eq:sbmbar2}, and let $\rho _0^*\left( \bar{\mathbf{x}}^*,\bar{\mathbf{y}}^*\right) $ be the optimal result of \eqref{eq:sbmbar}$'$, where $\left( \bar{\mathbf{x}}^*,\bar{\mathbf{y}}^*\right) =\left( \mathbf{x}+\mathbf{t}^{-*}, \mathbf{y}-\mathbf{t}^{+*}\right) $ is a super-efficiency projection of $\left( \mathbf{x},\mathbf{y}\right) $. On one hand, if $\mathbf{s}^{-*},\mathbf {s}^{+*}$ are optimal for \eqref{eq:sbmbar} then, by Proposition \ref{propn}, we have
\begin{equation}
\label{eq:32xx}
\rho ^*\left( \bar{\mathbf{x}}^*,\bar{\mathbf{y}}^*\right) =\rho \left( \bar{\mathbf{x}}^*,\bar{\mathbf{y}}^*, \mathbf{s}^{-*},\mathbf {s}^{+*}\right) =\rho \left( \mathbf{x}, \mathbf{y}, \mathbf{s}^{-*},\mathbf {s}^{+*}\right) \geq \rho _0^*\left( \bar{\mathbf{x}}^*,\bar{\mathbf{y}}^*\right) .
\end{equation}
On the other hand, if $\mathbf{s}^{-*},\mathbf {s}^{+*}$ are optimal for \eqref{eq:sbmbar}$'$ then, it can be proved analogously to Proposition \ref{propn} that if $t_i^{-*}>0$ for some $i\in \left\{1,\ldots ,m\right\}$, then $s_i^{-*}=0$, and if $t_r^{+*}>0$ for some $r\in \left\{1,\ldots ,s\right\}$, then $s_r^{+*}=0$. Hence
\begin{equation}
\label{eq:33xx}
\rho _0^*\left( \bar{\mathbf{x}}^*,\bar{\mathbf{y}}^*\right) =\rho \left( \mathbf{x}, \mathbf{y}, \mathbf{s}^{-*},\mathbf {s}^{+*}\right) =\rho \left( \bar{\mathbf{x}}^*,\bar{\mathbf{y}}^*, \mathbf{s}^{-*},\mathbf {s}^{+*}\right) \geq \rho ^*\left( \bar{\mathbf{x}}^*,\bar{\mathbf{y}}^*\right) .
\end{equation}
By \eqref{eq:32xx} and \eqref{eq:33xx} we have that $\rho ^*\left( \bar{\mathbf{x}}^*,\bar{\mathbf{y}}^*\right) =\rho _0^*\left( \bar{\mathbf{x}}^*,\bar{\mathbf{y}}^*\right) $. So, if $\mathbf{s}^{-*},\mathbf {s}^{+*}$ are optimal for \eqref{eq:sbmbar} we have that
$
\rho \left( \mathbf{x}, \mathbf{y}, \mathbf{s}^{-*},\mathbf {s}^{+*}\right) =\rho \left( \bar{\mathbf{x}}^*,\bar{\mathbf{y}}^*, \mathbf{s}^{-*},\mathbf {s}^{+*}\right) =\rho ^*\left( \bar{\mathbf{x}}^*,\bar{\mathbf{y}}^*\right) = \rho _0^*\left( \bar{\mathbf{x}}^*,\bar{\mathbf{y}}^*\right) 
$,
and then $\mathbf{s}^{-*},\mathbf {s}^{+*}$ are also optimal for \eqref{eq:sbmbar}$'$. Analogously, if $\mathbf{s}^{-*},\mathbf {s}^{+*}$ are optimal for \eqref{eq:sbmbar}$'$, then they are also optimal for \eqref{eq:sbmbar}.
\end{proof}

\begin{proof}[\textbf{Proof of Proposition \ref{prop41g2}}] Property 2. can be deduced directly from \eqref{eq:sbmcm}.

\begin{enumerate}
\item[1.] If $\left( \mathbf{x},\mathbf{y}\right) $ is inefficient, then $\delta ^*\left( \mathbf{x},\mathbf{y}\right) =1$. On the other hand, $\left( \mathbf{x},\mathbf{y}\right) \in \bar{P}$ and hence $\max \rho ^*|_{\bar{P}}=\rho ^*\left( \mathbf{x},\mathbf{y}\right)$ since $\rho ^*$ is strongly monotonic (although only weak monotonicity is needed) and any activity in $\bar{P}$ is dominated by $\left( \mathbf{x},\mathbf{y}\right) $. So, $\gamma \left( \mathbf{x},\mathbf{y}\right)  =\rho ^*\left( \mathbf{x},\mathbf{y}\right) $ by \eqref{eq:sbmcm}.


\item[3.] Since $\max \rho ^*|_{\bar{P}}<1$, by \eqref{eq:sbmcm}, $\gamma \left( \mathbf{x},\mathbf{y}\right) <\delta ^*\left( \mathbf{x},\mathbf{y}\right) $.
On the other hand, since $\left( \mathbf{x},\mathbf{y}\right) $ is efficient and there are not efficient activities in $\bar{P}$, we have that $\left( \mathbf{x},\mathbf{y}\right) \notin \bar{P}$ and hence, there are optimal super-efficiency slack vectors for \eqref{eq:ssbmact} that are not simultaneously zero. So, $\delta ^*\left( \mathbf{x},\mathbf{y}\right) >1$ and then $\gamma \left( \mathbf{x},\mathbf{y}\right) >\max \rho ^*|_{\bar{P}}$ by \eqref{eq:sbmcm}.
\end{enumerate}
\end{proof}

\begin{proof}[\textbf{Proof of Proposition \ref{propcont}}]
We have that $\max \rho ^*|_{\bar{P}}$ depends on $\left( \mathbf{x},\mathbf{y}\right) $ through $\bar{P}$ (i.e. the set of activities in $P$ that are dominated by $\left( \mathbf{x},\mathbf{y}\right) $). From \eqref{eq:barP}, it is clear that $\bar{P}$ varies in a continuous way with respect to $\left( \mathbf{x},\mathbf{y}\right) $. Moreover, since $\rho ^*|_P$ is continuous and $\bar{P}\subseteq P$, we have that $\max \rho ^*|_{\bar{P}}$ is continuous. Finally, since $\delta ^*$ is continuous, we conclude that $\gamma =\delta^* \cdot \max \rho ^*|_{\bar{P}}$ is also continuous.
\end{proof}

\begin{proof}[\textbf{Proof of Proposition \ref{prop42mon}}]
Let $\left( \mathbf{x},\mathbf{y}\right) $ be an activity strictly dominated by $\left( \mathbf{x}',\mathbf{y}'\right) $. Considering $\bar{P}$ and $\bar{P}'$ the sets of activities in $P$ that are dominated by $\left( \mathbf{x},\mathbf{y}\right) $ and $\left( \mathbf{x}',\mathbf{y}'\right) $ respectively, it holds that $\bar{P} \subseteq \bar{P}'$ and hence, $\max \rho ^*|_{\bar{P}}\leq \max \rho ^*|_{\bar{P}'}$. On the other hand, $\delta ^*\left( \mathbf{x},\mathbf{y}\right) \leq \delta ^*\left( \mathbf{x}',\mathbf{y}'\right) $ since $\delta ^*$ is weakly monotonic. So, $\gamma \left( \mathbf{x},\mathbf{y}\right) \leq \gamma \left( \mathbf{x}',\mathbf{y}'\right) $ by \eqref{eq:sbmcm}.
\end{proof}

\begin{proof}[\textbf{Proof of Proposition \ref{prop71}}]
We are going to prove that, in order to compute \eqref{eq:rhobarP2}, the objective function of \eqref{eq:rhobar} can be replaced by $\rho \left( \mathbf{x},\mathbf{y},\mathbf{s}^-,\mathbf{s}^+\right) $.
Let \eqref{eq:rhobar}$'$ be the program \eqref{eq:rhobar} with objective function $\rho \left( \mathbf{x},\mathbf{y},\mathbf{s}^-,\mathbf{s}^+\right) $, let $\rho _0^*\left( \mathbf{x}+\mathbf{t}^{-}, \mathbf{y}-\mathbf{t}^{+}\right) $ be the optimal result of \eqref{eq:rhobar}$'$, and
let $\bar{P}'$ be the subset of $\bar{P}$ whose activities are not strictly dominated by any other activity in $\bar{P}$. It is clear that if an activity $\left(\bar{\mathbf{x}}, \bar{\mathbf{y}}\right) \in \bar{P}$ maximizes $\rho ^*|_{\bar{P}}$ or $\rho _0^*|_{\bar{P}}$, then $\left(\bar{\mathbf{x}}, \bar{\mathbf{y}}\right) \in \bar{P}'$, because $\rho ^*|_{\bar{P}}$ and $\rho _0^*|_{\bar{P}}$ are strongly monotonic (see Proposition \ref{prop:rhomon}). Hence
\begin{equation}
\label{eq:31}
\max \rho ^*|_{\bar{P}}=\max \rho ^*|_{\bar{P}'},\qquad
\max \rho _0^*|_{\bar{P}}=\max \rho _0^*|_{\bar{P}'}.
\end{equation}
Let $\left( \bar{\mathbf{x}}, \bar{\mathbf{y}}\right) =\left( \mathbf{x}+\mathbf{t}^-,\mathbf{y}-\mathbf{t}^+\right) \in \bar{P}'$ and let $\left( \bm{\lambda}^*,\mathbf{s}^{-*},\mathbf{s}^{+*}\right) $ be an optimal solution for \eqref{eq:rhobar} or \eqref{eq:rhobar}$'$. We define $\bm{\tau}^-,\bm{\tau}^+$ such that $\tau_i^-=\max \left\{ 0,t_i^{-}-s_i^{-*}\right\} $, $\tau_r^+=\max \left\{ 0,t_r^{+}-s_r^{+*}\right\} $ for $i=1,\ldots ,m$ and $r=1,\ldots ,s$.
From the constraints of \eqref{eq:rhobar} or \eqref{eq:rhobar}$'$ we have
\begin{equation*}
\begin{array}{l}
\mathbf{x}+\bm{\tau}^-\geq \mathbf{x}+\left( \mathbf{t}^--\mathbf{s}^{-*}\right) =X\bm{\lambda}^*,\\
\mathbf{y}-\bm{\tau}^+\leq \mathbf{y}-\left( \mathbf{t}^+-\mathbf{s}^{+*}\right) =Y\bm{\lambda}^*,
\end{array}
\end{equation*}
and so $\left( \mathbf{x}+\bm{\tau}^-,\mathbf{y}-\bm{\tau}^+\right) \in P$.
Moreover, $\left( \mathbf{x}+\bm{\tau}^-,\mathbf{y}-\bm{\tau}^+\right) \in \bar{P}$ because $\tau _i^-,\tau _r^+\geq 0$.
Let us suppose that there exists $i\in \left\{ 1,\ldots ,m\right\} $ such that $s_i^{-*},t_i^->0$. Then $\tau _i^-<t_i^-$ and hence, taking into account that $\bm{\tau}^-\leq \mathbf{t}^-$ and $\bm{\tau}^+\leq \mathbf{t}^+$, we have that $\left( \mathbf{x}+\mathbf{t}^-,\mathbf{y}-\mathbf{t}^+\right) $ is strictly dominated by $\left( \mathbf{x}+\bm{\tau}^-,\mathbf{y}-\bm{\tau}^+\right) $, that is a contradiction because $\left( \mathbf{x}+\mathbf{t}^-,\mathbf{y}-\mathbf{t}^+\right) $ is not strictly dominated by any other activity in $\bar{P}$. Analogously for outputs. So, if $t_i^->0$ then $s_i^{-*}=0$, and if $t_r^+>0$ then $s_r^{+*}=0$, concluding that
\begin{equation}
\label{eq:33}
\rho \left( \bar{\mathbf{x}}, \bar{\mathbf{y}},\mathbf{s}^{-*},\mathbf{s}^{+*}\right) =\rho \left( \mathbf{x}+\mathbf{t}^-,\mathbf{y}-\mathbf{t}^+,\mathbf{s}^{-*},\mathbf{s}^{+*}\right) = \rho  \left( \mathbf{x},\mathbf{y},\mathbf{s}^{-*},\mathbf{s}^{+*}\right) .
\end{equation}
On one hand, if $\mathbf{s}^{-*},\mathbf{s}^{+*}$ are optimal for \eqref{eq:rhobar} then, by \eqref{eq:33} we have
\begin{equation}
\label{eq:34}
\rho ^*\left( \bar{\mathbf{x}}, \bar{\mathbf{y}}\right) = \rho \left(\bar{\mathbf{x}}, \bar{\mathbf{y}},\mathbf{s}^{-*},\mathbf{s}^{+*}\right)
 = \rho  \left( \mathbf{x},\mathbf{y},\mathbf{s}^{-*},\mathbf{s}^{+*}\right) 
\geq \rho _0^*\left( \bar{\mathbf{x}}, \bar{\mathbf{y}}\right) .
\end{equation}
On the other hand, if $\mathbf{s}^{-*},\mathbf{s}^{+*}$ are optimal for \eqref{eq:rhobar}$'$ then, by \eqref{eq:33} we have
\begin{equation}
\label{eq:35}
\rho _0^*\left( \bar{\mathbf{x}}, \bar{\mathbf{y}}\right) = \rho  \left( \mathbf{x},\mathbf{y},\mathbf{s}^{-*},\mathbf{s}^{+*}\right)
= \rho \left( \bar{\mathbf{x}}, \bar{\mathbf{y}},\mathbf{s}^{-*},\mathbf{s}^{+*}\right)
\geq \rho ^*\left( \bar{\mathbf{x}}, \bar{\mathbf{y}}\right) .
\end{equation}
Hence, by \eqref{eq:34} and \eqref{eq:35} we have that $\rho ^*\left( \bar{\mathbf{x}}, \bar{\mathbf{y}}\right) =\rho _0^*\left( \bar{\mathbf{x}}, \bar{\mathbf{y}}\right) $ for all $\left( \bar{\mathbf{x}}, \bar{\mathbf{y}}\right) \in \bar{P}'$, and then
$\rho ^*|_{\bar{P}'} = \rho _0^*|_{\bar{P}'}$.
Taking into account \eqref{eq:31}
we conclude that
$\max \rho ^*|_{\bar{P}}=\max \rho ^*|_{\bar{P}'}=\max \rho _0^*|_{\bar{P}'}=\max \rho _0^*|_{\bar{P}}$.
\end{proof}

\bibliography{references.bib}
\bibliographystyle{abbrvnat}

\end{document}